\documentclass[letter, 11 pt, reqno]{article}

\usepackage{xr}
    \usepackage{graphicx}
    \usepackage{styleset}
    \usepackage{macros}

\begin{document}

\title{Dimension Reduction of Generalized ASD Instantons}

\author{Dylan Galt\thanks{dylan.galt@stonybrook.edu} \;and Langte Ma\thanks{ltma@sjtu.edu.cn}}

\address{Stony Brook University, New York 11790, USA \\
Shanghai Jiao Tong University, Shanghai 200240, China
}



\maketitle

\begin{abstract}
We study generalized anti-self-dual instantons defined over Riemannian manifolds equipped with a parallel codimension-$4$ differential form. In particular, for product Riemannian manifolds possessing such a form, we study dimension reduction phenomena, finding a topological criterion for bundles which, when satisfied, allows for a complete characterization of dimension reduction for the corresponding moduli space of generalized ASD instantons. By establishing an integrability result for families of connections, we then deduce explicit descriptions for these moduli spaces, including those of Hermitian Yang--Mills connections, $G_2$-, and $\Spin(7)$-instantons. When one factor in the product is a $4$-manifold, we establish well-behaved compactifications for these moduli spaces. 
\end{abstract}

\maketitle

\tableofcontents

\section{Introduction}\label{s1}

\subsection{Background}

After the success of applying instanton theory to the study of topology and geometry in low dimensions, Donaldson--Thomas \cite{DT98} proposed a generalization to higher-dimensional manifolds equipped with special geometric structures. As explored by Tian \cite{T00} and Tao--Tian \cite{TT04}, the major issue of extracting geometric information from the moduli space of instantons over these higher-dimensional manifolds arises from their complicated ‘bubbling sets,' which in general are merely rectifiable. Later, Donaldson--Segal \cite{DS11} elaborated on the initial program for this higher-dimensional instanton theory by studying the calibration property of bubbling sets, in the sense of \cite{HL82}. This has led to fruitful progress on the interplay between instantons and calibrated submanifolds. 

In this paper, we are mainly concerned with dimension reduction phenomena in instanton theory. This means that given a fibration, we seek to describe instantons—absolute minimizers of the Yang-Mills functional—over the total space of the fibration in terms of instantons over its fibres. Over the past few decades, study of this problem has generally focused on the case when the fibre has codimension one, exploiting the fact that in this case the defining equation for instantons is gauge-equivalent to the gradient flow equation of some functional on the space of connections on the fibres. This principle has been established rigorously in various gauge theory contexts, including for Seiberg--Witten monopoles on $3$-manifolds \cite{MOY97} and $4$-manifolds \cite{B01}, instanton Floer homology \cite{M99}, anti-self-dual (ASD) instantons on $4$-manifolds \cite{A88}, and $G_2$- and $\Spin(7)$-instantons \cite{W20}.

When it comes to fibrations with fibres of higher codimension, establishing a complete picture of dimension reduction becomes more delicate (see for instance \cite{W22}). However, it is quite common to encounter interesting fibrations of this type in the study of higher dimensional gauge theory, for instance with the resolution block of Joyce’s \cite{J96} compact $G_2$- and $\Spin(7)$-manifolds, and tubular neighborhoods of bubbling sets \cite{W17, W17a}. Thus, it is natural to seek an approach to studying dimension reduction in the absence of the gradient flow picture.

In this paper, we develop a theory for dimension reduction beyond codimension-1, and consequently establish global descriptions of the moduli space of instantons over several different families of product manifolds with calibrated factors. In most cases, we provide the first known complete description of the corresponding moduli space of instantons and our results essentially characterize all $G_2$- and $\textup{Spin}(7)$-instantons over smooth, compact, hyperkähler coassociative fibrations. Our approach is fundamentally geometric and particularly suited for future applications to adiabatic limit problems.

\subsection{Main Results}
The data of our general setup consists of a closed Riemannian manifold $M$ of dimension $n \geq 4$, and a differential  $(n-4)$-form $\theta \in \Omega^{n-4}(M)$. We consider the following pointwise defined, bounded, self-adjoint operator 
\begin{equation}
\begin{split}
\mathcal{Q}_{\theta}: \Omega^2(M) & \longrightarrow \Omega^2(M) \\
\omega & \longmapsto \star(\theta \wedge \omega),
\end{split}
\end{equation}
which has a bounded and closed spectral set $\Spec(\mathcal{Q}_{\theta}) \subseteq \R$. We denote the minimum of the spectral set by 
\begin{equation}
\lambda_{\theta} = \min \Spec(\mathcal{Q}_{\theta}). 
\end{equation}

Let $\pi: P \to M$ be a principal bundle with structure group $U(r)$ or $SU(r)$, for $r \in \Z_+$. We shall study connections $A$ on $P$ satisfying 
\begin{equation}\label{e1.3}
\star \left(\theta \wedge F^o_A \right) = \lambda_{\theta} \cdot F^o_A, 
\end{equation}
where $F^o_A:= F_A - (1/r) \tr(F_A) \cdot \id$ denotes the traceless part of the curvature $2$-form. For convenience, we shall refer to connections satisfying \eqref{e1.3} as {\em $\theta$-instantons} (we adopt this terminology from Tian \cite{T00}, who referred to solutions to \eqref{e1.3} as $\theta$-ASD-instantons). Denote by $\G_P:= \Aut(P)$ the gauge group, which acts on connections by pull-back. The moduli space of $\theta$-instantons is defined to be 
\begin{equation}
\M_{\theta}(M, P):= \{\text{connections on } P \text{ satisfying } \eqref{e1.3}\}/ \G_P.
\end{equation}

\begin{dfn}\label{d1}
We say that the $(n-4)$-form $\theta$ on $M$ is {\em admissible} if the following two conditions hold:
\begin{enumerate}
\item $\theta$ is parallel with respect to the Levi--Civita connection on $M$. 
\item The minimum of the spectrum $\lambda_{\theta} = -1$. 
\end{enumerate}
\end{dfn}

These conditions ensure that $\lambda_{\theta}$ is an eigenvalue of $\mathcal{Q}_{\theta}$ and that $\theta$-instantons minimize the Yang--Mills functional (see \autoref{s2}). Note that (b) is merely a normalization requirement, which can always be arranged by rescaling $\theta$. In contrast, condition (a) is quite restrictive, but is satisfied by all of the natural examples motivating our study. 

We now restrict ourselves to the case when $M$ is a product (which we regard as a trivial fibre bundle) and $\theta$ is adapted in a suitable sense to this product splitting. More precisely, let $M = Y \times X$ be the product of two closed Riemannian manifolds with $\dim X \geq 4$, and assume that the admissible $(n-4)$-form $\theta$ takes the form
\begin{equation}
\theta = \vol_Y \wedge \beta + \alpha.
\end{equation}
where $\beta$ is a codimension-$4$ differential form on $X$, and $\alpha$ is a differential form with mixed terms that involve those of lower degree along the $Y$-factor. 

It is not hard to show that $\lambda_{\beta} \geq \lambda_{\theta}$ (see \autoref{s4}) and so we may normalize $\theta$ by multiplying a real factor to guarantee that 
\begin{equation}
0 > \lambda_{\beta} \geq \lambda_{\theta} = -1. 
\end{equation}
When $\lambda_{\beta} = \lambda_{\theta} = -1$, we say that $\theta$ is {\em compatible} with the product structure of $M$. A special case of this is when $\alpha=0$, in which case $Q_\theta(\omega)=\lambda\omega$ for $\lambda\neq 0$ if and only if $\omega$ is pulled back from $X$ and $Q_\beta(\omega)=\lambda \omega$, and so every $\theta$-instanton on $M$ is the pullback of a $\beta$-instanton on $X$. Thus, the real difficulty in approaching the issue of dimension reduction at this level of generality is to understand which conditions must be placed on $\alpha$.

The cornerstone of our dimension reduction theory is a particular energy condition that intertwines the topology of the bundle $P$ and the cohomology class of the form $\alpha$. 

\begin{dfn}
We define the {\em $\alpha$-charge} of the bundle $P$ to be the pairing 
\begin{equation}
\kappa_{\alpha}(P):= \langle \kappa(P) \smile [\alpha], [M] \rangle \in \R,
\end{equation}
where $\kappa(P) = c_2(P) - (1/2) c_1^2(P) \in H^4(M; \Z)$. Depending on $\kappa_{\alpha} < 0$, $\kappa_{\alpha} = 0$, or $\kappa_{\alpha} > 0$, we say $P$ is $\alpha$-negative, $\alpha$-null, or $\alpha$-positive, respectively.    
\end{dfn}

The geometric relevance of $\alpha$-charge to the $\theta$-instanton equation can be seen from the following perspective. We write $\pi_X: P_X \to X$ as a reference bundle isomorphic to the restriction of $P$ to the fibres $X_y:=\{y\}\times X$. Denote by $\A_{P_X}$ the space of connections on $P_X$ and $\B_{P_X} = \A_{P_X} / \Aut(P_X)$ the configuration space. By restricting to the fibres $X_y$, a connection $A$ on $P$ defines a tautological map
\begin{equation}
\mathscr{F}_A: Y \longrightarrow \B_{P_X}. 
\end{equation}
The Dirichlet energy
\begin{equation}
E(\mathscr{F}) := \frac{1}{2} \int_Y |d\mathscr{F}|^2 \vol_Y
\end{equation}
then defines a functional on the space of smooth maps $\Map(Y, \B_{P_X})$ and, in favorable cases (see \cite{GM24}), one can show that this attains its absolute minimum at so-called `adiabatic limit $\theta$-instantons' with value given by precisely the $\alpha$-charge $\kappa_{\alpha}(P)$. 

Thus, the picture of dimension reduction we seek to establish in this paper is intimately connected with the idea of adiabatic limits and can be summarized, roughly speaking, as follows: one can imagine collapsing the moduli space of $\theta$-instantons over $M$ down to the moduli space of $\beta$-instantons over $X$, and the `energy lost' in this process is encoded by the $\alpha$-charge of $P$. If $P$ is $\alpha$-null, there is essentially no energy lost and one expects a complete picture of dimension reduction, fully recovering the moduli spaces of $\theta$-instantons over $M$ from the moduli space of $\beta$-instantons over $X$. Our main results establish this rigorously.

To simplify the exposition, we work with $SU(r)$-bundles and irreducible instantons first. We say an $SU(r)$-connection $A$ on $P$ is {\em irreducible} if its stabilizer group $\Stab(A)$ is the center $\Z/r$ of $SU(r)$. We write $\M^*_{\theta}(M, P)$ for the moduli space of irreducible $\theta$-instantons. The following is our main dimension reduction theorem.

\begin{thm}\label{t1.1}
Let $M = Y \times X$ be the product of two closed Riemannian manifolds and $\theta = \vol_Y \wedge \beta + \alpha$ an admissible $(n-4)$-form compatible with the product structure on $M$. Suppose $\pi: P \to M$ is an $\alpha$-null principal $SU(r)$-bundle such that $\M^*_{\theta}(M, P) \neq \varnothing$. Then there is a canonical isomorphism of irreducible moduli spaces 
\[
\Theta^*: H^1(Y; \Z/r) \times \mathcal{N}^*_{\beta}(X, P_X) \longrightarrow \mathcal{M}^*_{\theta}(M, P),
\]
where $\mathcal{N}^*_{\beta}(X, P_X)$ denotes the moduli space of irreducible $\beta$-instantons over $P_X$. 
\end{thm}

In \autoref{s4}, we shall give a necessary and sufficient criterion for the non-vanishing of moduli spsaces over $\alpha$-null bundles. In the case of pull-back bundles, it is not hard to see that $\M^*_{\theta}(M, P) \neq \varnothing$ if and only if $\mathcal{N}^*_{\beta}(X, P_X) \neq \varnothing$. 

When $P$ is $\alpha$-negative, or $\alpha$-null but $\lambda_{\beta} > \lambda_{\theta}$, we prove in \autoref{s4} that the moduli space $\M_{\theta}(M, P)$ is empty. Due to elliptic regularity, we are considering $C^{\infty}$-connections and so the moduli spaces are equipped with the $C^{\infty}$-topology. The map $\Theta^*$ is an isomorphism between analytic varieties, which means that it is a homeomorphism and identifies Zariski tangent spaces. The first factor $H^1(Y; \Z/r)$ corresponds to flat connections on the trivial line bundle over $Y$ with holonomy contained in $\Z/r$. The map $\Theta^*$ is constructed in \autoref{s4} by twisting connections on $P$ by the pull-back of such flat connections.

When $r=2$, the other stabilizer that a non-trivial connection can carry is the maximal torus $U(1) \subseteq SU(2)$. We shall refer to these connections as {\em abelian} connections. We write a superscript $^a$ in $\M^a_{\theta}$ and $\mathcal{N}^a_{\beta}$ to indicate the moduli spaces consisting of abelian instantons. The following is our main dimension reduction theorem for these reducible connections.

\begin{thm}\label{t1.2}
Let $(M, \theta)$ be a pair as in \autoref{t1.1} and $\pi: P \to M$ an $\alpha$-null principal $SU(2)$-bundle such that $\M^a_{\theta}(M, P) \neq \varnothing$. Then there is a canonical isomorphism between the abelian loci
\[
\Theta^a: \Pic(Y) \times \mathcal{N}^a_{\beta}(X, P_X) \longrightarrow \mathcal{M}^a_{\theta}(M, P),
\]
where $\Pic(Y)$ is the Picard torus of $Y$ that parametrizes flat connections on the trivial line bundle over $Y$.
\end{thm}

In principle, given any fixed subgroup $\Gamma \subseteq SU(r)$ that arises as the stabilizer of some connection, we can define the twisting map $\Theta^{\Gamma}$ and identify the strata of the moduli spaces consisting of connections whose stabilizers are $\Gamma$. However, it is a cumbersome task to exhaust all possibilities. Since we do not have in mind any further applications for such refinements, we shall be satisfied with the current description rather than an explicit stratification in general. 

By allowing for connections with rectifiable singular loci, Tian \cite{T00} extended Uhlenbeck's compactness result \cite{U82} to $\theta$-instantons over higher-dimensional manifolds. In the case of $\alpha$-null bundles $P \to Y \times X$ with $\dim X = 4$, we further improve Tian's compactification by proving that the singular loci are well-behaved and can be removed in the limiting instantons. 

We assume that $X$ is a closed $4$-manifold and $\beta \equiv 1$ in the pair $(M, \theta)$. Given an $\alpha$-null principal $SU(r)$-bundle $P \to M$ and an integer $1 \leq s \leq \kappa(P_X)$, we shall construct another $\alpha$-null bundle $P(s)$ in \autoref{s4} characterized by the change on its second Chern class
\begin{equation}
c_2(P) - c_2(P(s)) = s\PD([Y]) \in H^4(M; \Z). 
\end{equation}
An {\em ideal $\theta$-instanton} on the $\alpha$-null bundle $P$ is a pair $(A_o, \pmb{x})$ consisting of a $\theta$-instanton $A_o$ on $P(s)$ and an element $\pmb{x} \in \Sym^s(X)$, where $\Sym^s(X)$ is the $s$-fold symmetric product of $X$. 

\begin{dfn}\label{d1.3}
We say a sequence of $\theta$-instantons $\{A_n\}_{n\geq 1}$ on an $\alpha$-null bundle $P$ converges to an ideal $\theta$-instanton $(A_o, \pmb{x})$ if 
\begin{enumerate}
\item for any $C^{\infty}$ function $f: M \to \R$, one has 
\[
\int_{M} f \cdot |F_{A_n}|^2 \vol_M \longrightarrow \int_{M} f \cdot |F_{A_o}|^2 \vol_M + \sum_{x \in \pmb{x}} \int_{Y \times \{x\}} f \vol_Y \text{ as $n \to \infty$}. 
\]
\item one can find bundle isomorphisms $v_n: P(s)|_{M \backslash Y \times \pmb{x}}  \to P|_{M \backslash Y \times \pmb{x}}$ so that 
\[
v^*_nA_n \longrightarrow A_o \text{ in $C^{\infty}(K)$ as $n \to \infty$}
\]
for any compact subset $K \subseteq M \backslash Y \times \pmb{x}$.
\end{enumerate}
\end{dfn}

We establish the following compactification result.

\begin{thm}\label{t1.4}
Let $(M, \theta)$ be pair in \autoref{t1.1} and assume further that $X$ is a closed $4$-manifold and $\beta \equiv 1$. Given an $\alpha$-null principal $SU(r)$-bundle, the space 
\[
\overline{\M}_{\theta}(M, P) :=  \M_{\theta}(M, P)\cup \bigcup_{s=1}^{\kappa(P_X)} \M_{\theta}(M, P(s)) \times \Sym^s(X). 
\]
of ideal $\theta$-instantons on $P$ is compact. 
\end{thm}

Note that Tian \cite{T00} proved that the bubbling locus for $\theta$-instantons must be calibrated by $\theta$ and in fact \autoref{t1.4} tells us that the bubbling locus for instantons on $\alpha$-null bundles must be the slices $Y \times \{x\}$, not any other calibrated submanifolds in $M$. Moreover, appealing to the gluing theorem of Taubes \cite{T82}, bubbling along any slice $Y \times \{x\}$ can be realized over such $\alpha$-null bundles, in contrast with Walpuski’s result \cite{W17, W17a}. This is because our form $\theta$ is not `generic' and for index reasons our instantons are obstructed.

We will apply our main dimension reduction results to the following concrete examples of manifolds, most of which have special or exceptional holonomy:
\begin{enumerate}
\item flat tori $T^n$ of dimension $n$;
\item K\"ahler manifolds $Z$ of dimension $n = 2m$ equipped with a K\"ahler form $\omega \in \Omega^{(1,1)}_{\R}$;
\item hyperk\"ahler $4$-manifolds $X$ equipped with a triple of K\"ahler forms $(\omega_1, \omega_2, \omega_3)$ satisfying 
\[
\omega_i \wedge \omega_j = 2\delta_{ij} \vol_X;
\]
\item Calabi--Yau $3$-folds $Z$ equipped with a K\"ahler form $\omega \in \Omega^{(1, 1)}_{\R}$ and a holomorphic volume form $\Omega \in \Omega^{(3,0)}$;
\item $G_2$-manifolds $V$ of dimension $7$ equipped with an associative calibration form $\varphi \in \Omega^3$ and its Hodge dual $\psi = \star \varphi \in \Omega^4$;
\item $\Spin(7)$-manifolds $M$ of dimension $8$ equipped with a self-dual Cayley calibration form $\Phi \in \Omega^4$. 
\end{enumerate}

Here, we have put parentheses in the superscript of differential forms to indicate that the bigrading is induced from the complex structure. In \autoref{t1} below, we list the types of factors in the product manifold $M = Y \times X$ with which we are concerned. 

\begin{table}[h]
\begin{tabular}{| c | c | c | c | c | }
 \hline
 $M$ & $Y$ & $X$ & $\theta$ & $\alpha$ \\ 
 \hline
$G_2$ & $T^3$ & $\omega_1, \omega_2, \omega_3$ & $\vol_Y + \sum_i dt_i \wedge \omega_i$ & $\sum_i dt_i \wedge \omega_i$ \\
 \hline
 $G_2$ & $S^1$ & $\omega, \Omega$ & $\vol_Y \wedge \omega + \Rea \Omega$ & $ \Rea \Omega$ \\
 \hline
  $\Spin(7)$ & $\tau_1, \tau_2, \tau_3$ &  $\omega_1, \omega_2, \omega_3$ & $\vol_Y + \vol_X - \sum_i \tau_i\wedge \omega_i$ & $\vol_X - \sum_i \tau_i\wedge \omega_i$   \\
 \hline
   $\Spin(7)$ & $S^1$ & $\varphi, \psi$ & $\vol_Y \wedge \varphi + \psi$ & $\psi$ \\
 \hline
K\"ahler & $\omega_1$ & $\omega_2$ & $(\omega_1 + \omega_2)^{m-2}$ & $\sum_{k \neq m_1} c_k \omega_1^k \wedge \omega_2^{m - 2 - k}$ \\
\hline
\end{tabular}
\vspace{3mm}
\caption{Product Factors}
\label{t1}
\end{table}

In the first row of \autoref{t1}, the hyperk\"ahler fibres $X_y$ are calibrated by the form $\star \theta$ and are called coassociative. In the third row, the hyperk\"ahler fibres $X_y$ are calibrated by $\theta$ and are called Cayley. In the last row of \autoref{t1}, $Y$ and $X$ are K\"ahler manifolds of dimension $m_1$ and $m_2$ respectively, $m = m_1 + m_2$, and $c_k = {m-2 \choose k}$. The last case provides us with examples where the codimension of the dimension reduction can be arbitrarily large.

It is evident that pull-back bundles $\pi_X^* P_X \to M$ are $\alpha$-null bundles. Of course there are lots of $\alpha$-null bundles that are not topologically isomorphic to a pull-back bundle whose characteristic classes have to be pulled-back ones. Nevertheless, if we work with low-dimensional $Y$-factor, pull-back bundles are necessary to guarantee the existence of $\theta$-instantons among all $\alpha$-null bundles. 

\begin{prop}\label{p1.5}
Let $(M, \theta)$ be a pair as before with $\dim Y \leq 4$. Suppose $P \to M$ is an $\alpha$-null principal $SU(r)$-bundle satisfying $\M_{\theta}(M, P) \neq \varnothing$. Then $P$ is topologically isomorphic to a pull-back bundle $\pi_2^*P_X$ with respect to some $SU(r)$-bundle $P_X \to X$. 	
\end{prop}

The strategy for the proof of our main theorem is to combine the $\alpha$-null energy condition with an integrability argument for families of connections. To formulate this, we consider an arbitrary principal $G$-bundle $\pi: P \to M = Y \times X$ with $G$ a compact linear algebraic group. We obtain a natural fibre bundle $\underline{\pi}: \underline{P} \to Y$ with fibres $P|_{X_y}$, where $\underline{\pi}=\pi_Y\circ \pi$, which we refer to as the corresponding \textit{family bundle}. The structure group of the family bundle is the gauge group $\G_{P_X} = \Aut(P_X)$. We denote by $\Stab(B) \subseteq \G_{P_X}$ the stabilizer group of a connection $B$ on $P_X$ and for convenience we refer to principal $G$-connections on $P$ as $G$-connections.

The product structure of $M$ decomposes differential forms into a bi-graded sum 
\begin{equation}\label{e1.1}
\Omega^l(M) = \bigoplus_{s+t = l} \Omega^{s, t}(M). 
\end{equation}
In particular, the curvature form $F_A$ of a $G$-connection $A$ decomposes as
\[
F_A = F^{2,0}_A + F^{1,1}_A + F^{0,2}_A. 
\]
Fixing a basepoint $y_0 \in Y$, we say a gauge transformation $u \in \G_{P}$ is \textit{fibre-based} at $y_0 \in Y$ if $u|_{P_{y_0}} = \id$. The following is our main integrability result.

\begin{thm}\label{t1.6}
Let $P \to Y \times X$ be a principal $G$-bundle and $A$ a smooth $G$-connection on $P$ such that 
\[
F^{2, 0}_A = 0 \qquad F^{1, 1}_A = 0. 
\]
Then one can find a representation $\rho: \pi_1(Y, y_0) \to \Stab(A|_{P_{y_0}})$ satisfying:
\begin{enumerate}
\item[\upshape{(a)}] The family bundle $\underline{P}$ is isomorphic to $\tilde{Y} \times_{\rho} P_X$ where $\tilde{Y} \to Y$ is the universal cover of $Y$.
\item[\upshape{(b)}] The connection $A$ is fibre-based gauge-equivalent to the connection obtained by the parallel transport of $A|_{P_{y_0}}$ using the flat connection on $\underline{P}$ specified by $\rho$. 
\end{enumerate}
\end{thm}

Our main results expand upon several theorems in the literature. In \cite{W20}, Wang studied the codimension-$1$ reduction of $G_2$- and $\Spin(7)$-instantons in the case of pull-back $U(r)$-bundles. His formulation is different from ours, exploiting the Chern--Simons functional and introducing `broken’ gauge transformations to avoid the discussion of monodromy arising from the base $Y=S^1$. In the setting of codimension-1 reduction, the twisting map $\Theta$ clarifies the role of this monodromy and identifies the Zariski tangent spaces of both moduli spaces. The approach adopted in this paper works for bundles with more general structure groups as well. In \cite{S14}, Sá Earp studied $G_2$-instantons over associative torus fibrations. Interestingly, \cite[Theorem 2]{S14} provides us with non-standard $G_2$-structures on $T^7 = T^3 \times T^4$ for which pull-back bundles from $T^4$ admit no $G_2$-instantons. This is evidence that the moduli space of $\theta$-instantons is obstructed and is sensitive to the fibration structure on $(M, \theta)$.

This paper leaves open the issue of characterizing the moduli space of $\theta$-instantons over $\alpha$-positive bundles. Consistent with the aforementioned adiabatic picture, when the $\alpha$-charge $\kappa_{\alpha}$ is sufficiently small we expect that solving this problem involves understanding whether one can reconstruct $\theta$-instantons from adiabatic limit $\theta$-instantons. We emphasize this perspective in a forthcoming second paper \cite{GM24}. However, it is still mysterious to us what phenomena will occur as $\kappa_{\alpha}$ gets large. On the other hand, we believe that \autoref{t1.1}, suitably phrased, should hold for fibrations with non-compact fibres. If one restricted oneself to instantons with finite energy, the generalization would not be very difficult and could be applied to Walpuski’s construction \cite{W13b} of instantons over generalized Kummer constructions. However, it is interesting to work out more delicate control on the energy growth along the non-compact ends, which matches better with the Floer theory picture on cylinders. We plan to take up this issue in a future work.

\subsection*{Acknowledgments}
Both authors are grateful to Simon Donaldson for his encouragement and enlightening discussions. This project began when the second author held a postdoctoral position at the Simons Center for Geometry and Physics. He would like to acknowledge the Simons Center with gratitude for its support and hospitality. This bulk of this project was carried out while the first author was a collaboration graduate student in the Simons Collaboration on Special Holonomy in Geometry, Analysis, and Physics. He would like to thank the collaboration for its generous support.

\section{Generalized Anti-Self-Dual Instantons}\label{s2}

Generalized anti-self-dual instantons were introduced by Tian \cite{T00} following the work of Donaldson--Thomas \cite{DT98}. In this section, we shall set up carefully the version of this framework suited to this paper, with an emphasis on the Yang--Mill functional and local deformations. We will also review examples arising from special holonomy geometry, which initially motivated this study.

\subsection{Anti-Self-Dual $\theta$-Instantons}

Let $M$ be a smooth, oriented, Riemannian $n$-manifold with $n \geq 4$, and $\theta \in \Omega^{n-4}(M)$ a non-zero $(n-4)$-form. We consider the following operator on the space of $2$-forms
\begin{equation}
\begin{split}
\mathcal{Q}_{\theta}: \Omega^2(M) & \longrightarrow \Omega^2(M) \\
\omega & \longmapsto \star(\theta \wedge \omega). 
\end{split}
\end{equation}
It is easy to see that $\mathcal{Q}_{\theta}$ is self-adjoint and bounded with respect to the $L^2$-inner product on $\Omega^2(M)$ induced by the Riemannian structure on $M$, and thus admits a real, bounded, and closed spectrum $\Spec(\mathcal{Q}_{\theta}) \subseteq \R$.

\begin{lem}\label{l2}
If $\theta$ is non-zero, the operator $\mathcal{Q}_{\theta}$ is not identically zero.
\end{lem}
\begin{proof}
It is enough to construct $\omega,\omega'\in \Omega^2(M)$ such that $\int_{M}\theta\wedge\omega\wedge\omega'\neq 0$. Let $p\in M$ be such that $\theta_p\neq 0$ and let $U$ denote a coordinate chart centered at $p$ with coordinates $x_1,\cdots, x_n$. We can write $\theta|_U = \sum_{I}a_Idx_I$ where $I$ denotes a multi-index satisfying $|I|=n-4$. Since $\theta_p\neq 0$, there is some $I$ such that $a_I(p)\neq 0$, and shrinking $U$ if necessary we can arrange that $a_I|_U$ is non-vanishing. Let $J$ denote the multi-index dual to $I$ and choose some $\omega_0,\omega'_0\in\Omega^2(U)$ such that $\omega_0\wedge\omega'_0=dx_J$. Let $f: M\rightarrow\mathbb{R}$ be a bump function such that $f|_{U'}\equiv 1$ on an open set $U'\subset U$ and $f|_{M\setminus U}\equiv 0$. Then, $\omega=f\cdot \omega_0$ and $\omega'=f\cdot\omega'_0$ are well-defined on all of $M$ and by construction $\int_M\theta\wedge\omega\wedge\omega'=\int_U f^2a_I dx_1\wedge\cdots\wedge dx_n\neq 0$.
\end{proof}

It follows from \autoref{l2} that the minimum of $\Spec(\mathcal{Q}_{\theta})$ is nonzero and so by multiplying $\theta$ by a real factor we can arrange that it is $-1$. We denote it by
\begin{equation}
\lambda_{\theta} = \min \Spec(\mathcal{Q}_{\theta}) = -1. 
\end{equation}

To simplify the matter, we further assume that $\lambda_{\theta}$ is an eigenvalue of $\mathcal{Q}_{\theta}$, i.e. $\ker (\mathcal{Q}_{\theta} - \lambda_{\theta} \id) \neq \{0\}$. The simplest way to achieve this is to assume that $\theta$ is parallel, hence admissible in the sense of \autoref{d1}. When $\theta$ is parallel, the bundle $\Lambda^2T^*M$ of alternating $2$-tensors can be decomposed into sub-bundles corresponding to the eigenspace decomposition of the $\mathcal{Q}_{\theta}$-action on $\Lambda^2T^*_pM$ at a single point $p \in M$. Denoting the eigenvalues of $\mathcal{Q}_{\theta}$ by $\lambda_{\theta} = \lambda_0 < \lambda_1 < ... < \lambda_k$, we have the corresponding sub-bundles 
\begin{equation}
\Lambda^2_{\lambda_i}T^*M:= \{ \omega \in \Lambda^2T^*M: \star(\theta \wedge \omega) = \lambda_i \omega \}.
\end{equation}
We shall denote orthogonal projection onto these sub-bundles by
\begin{equation}
\pi_{\lambda_i}: \Lambda^2T^*M \longrightarrow \Lambda^2_{\lambda_i}T^*M. 
\end{equation}
Writing $\pi_{\theta} = \pi_{\lambda_0}$ for the projection onto the $-1$-eigen-sub-bundle, and $\vol_M$ is the volume form on $M$, we have the following useful identity
\begin{equation}\label{e2.5}
\omega \wedge \omega \wedge \theta = \left(\sum_{i=1}^k \lambda_i |\pi_{\lambda_i}(\omega)|^2 - |\pi_{\theta}(\omega)|^2 \right) \vol_M.
\end{equation} 

Henceforth, unless otherwise stated, we assume that $\theta$ is admissible and that $M$ is closed, referring to such a pair $(M,\theta)$ as an \textit{admissible pair}. Switching to gauge theory, we apply the previous discussion to Lie algebra-valued forms. Let $\pi: P \to M$ be a principal $SU(r)$-bundle and denote by $\A_P$ the space of all smooth $SU(r)$-connections on $P$. We write $\g_P:=P \times_{\ad} \mathfrak{su}(r)$ for the induced adjoint bundle and note that $\A_P$ is an affine space modeled on $\Omega^1(M, \g_P)$, which in turn is equipped with a natural metric by coming from the Killing form on $\mathfrak{su}(r)$: 
\begin{equation}
\langle a, b \rangle :=-\tr (\ad(a) \circ \ad(b)) = -2r\tr(ab), \quad \forall a, b\in \su(r). 
\end{equation}

We denote by $F_A \in \Omega^2(M, \g_P)$ the curvature of an $SU(r)$-connection $A$ on $P$. The Yang--Mills functional on the connection space $\A_P$ is then defined as 
\begin{equation}
\mathcal{YM}(A):= \int_M  |F_A|^2 \vol_M.
\end{equation}
The critical points of this functional (cf. \cite{DK90}) are those connections satisfying the Yang--Mills equations
\begin{equation}
d^*_A F_A = 0. 
\end{equation}
If the curvature $F_A \in \Omega^2(M, \g_P)$ is an eigenvector of the operator $\mathcal{Q}_{\theta}$, then
\begin{equation}
\lambda_i \cdot d^*_A F_A = d^*_A \star(\theta \wedge F_A) = 0, 
\end{equation}
using the Bianchi identity and the fact that $\theta$ is parallel, hence in particular closed. If $\lambda_i \neq 0$, then we can conclude that $A$ is a Yang--Mills connection. 

\begin{dfn}
An $SU(r)$-connection $A$ on $P$ is called an (anti-self-dual) $\theta$-instanton if 
\[
\star(\theta \wedge F_A) = -F_A. 
\]
\end{dfn}

It is well-known that anti-self-dual instantons over $4$-manifolds absolutely minimize the Yang--Mills functional and that any such absolute minimizer must be anti-self-dual. Appealing to \eqref{e2.5}, we see that $\theta$-instantons retain this characterizing feature. Indeed, this is the primary reason we ask $F_A$ to be an eigenvector for $Q_\theta$ corresponding to the minimal eigenvalue $\lambda_\theta=-1$.

\begin{prop}\label{p2.3}
Let $(M, \theta)$ be an admissible pair, and $\pi: P \to M$ a principal $SU(r)$-bundle. Then any $\theta$-instanton $A$ minimizes the Yang--Mills functional on $\A_P$. Moreover, the minimum value of the Yang-Mills functional is given by the topological quantity
\[
\kappa_{\theta}(P):= 16r\pi^2 \langle c_2(P) \smile [\theta], [M] \rangle \in \R. 
\]
\end{prop}

\begin{proof}
Since any parallel form is harmonic, we know that $\theta$ is closed. Thus $\kappa_{\theta}(P)$ is a topological quantity that depends only on the cohomology class of $\theta$. 

Let $A$ be an $SU(r)$-connection on $P$. Then \eqref{e2.5} implies that 
\[
-2r\tr(F_A \wedge F_A) \wedge \theta = \left( \sum_{i=1}^k \lambda_i \cdot |\pi_{\lambda_i}(F_A)|^2 - |\pi_{\theta}(F_A)|^2 \right) \vol_M. 
\]
The Yang--Mills functional can be computed as 
\[
\begin{split}
\mathcal{YM}(A) & = \sum_{i=1}^k \| \pi_{\lambda_i}(F_A) \|^2_{L^2} + \| \pi_{\theta}(F_A) \|^2_{L^2} \\
& = \sum_{i=1}^k (1 + \lambda_i) \cdot \| \pi_{\lambda_i}(F_A) \|^2_{L^2} + 2r\int_M \tr(F_A \wedge F_A) \wedge \theta \\
& \geq 16r\pi^2 \langle c_2(P) \smile [\theta], [M] \rangle,
\end{split}
\]
where the last inequality makes use of the fact that $\lambda_i > -1$ for $i =1, ..., k$. It is clear that equality holds if and only if that $\pi_{\theta}(F_A) = F_A$, which means that $A$ is a $\theta$-instanton. 
\end{proof}

\begin{rem}
The admissibility assumption on $\theta$ in \autoref{p2.3} can certainly be weakened. All we need is that $\theta$ is closed and that the minimum $\lambda_{\theta}$ of $\textup{Spec}(\mathcal{Q}_\theta)$ is an eigenvalue of $\mathcal{Q}_{\theta}$.
\end{rem}

We continue on to discuss the deformation theory of $\theta$-instantons. The gauge group $\G_P:= \Aut(P)$ acts on the space of connections by pulling back
\begin{equation}
u \cdot A:= A - (d_A u)u^{-1}, \quad u \in \G_P, \;\; A \in \A_P, 
\end{equation}
and it is straightforward to verify that the defining equation for $\theta$-instantons is $\G_P$-equivariant. Thus we define the moduli space of $\theta$-instantons to be 
\begin{equation}
\M_{\theta}(M, P):= \{ A \in \A_P: \star(\theta \wedge F_A) = -F_A\}/\G_P. 
\end{equation}

After linearizing the gauge group action and the defining equation, we get the following differential complex coupled to the $\theta$-instanton $A$
\begin{equation}\tag{$E_{\theta, A}(M)$}
0 \to \Omega^0(M, \g_P) \xrightarrow{-d_A} \Omega^1(M, \g_P) \xrightarrow{\pi_{\theta^{\perp}} \circ \;d_A} \Omega^2_{\theta^{\perp}}(M, \g_P) \to 0,
\end{equation}
where $\Omega^2_{\theta^{\perp}}(M, \g_P):=\Omega^0(M,\ker(\pi_\theta))$. Denote $L^2$-orthogonal projection onto this subspace $\Omega^2_{\theta^{\perp}}(M, \g_P)$ by $\pi_{\theta^\perp}$. We shall refer to $E_{\theta, A}(M)$ as the deformation complex of the moduli space $\M_{\theta}(M, P)$ at $[A]$. In general, $E_{\theta, A}(M)$ is not an elliptic complex, but we shall see later that for many important, concrete examples $E_{\theta, A}(M)$ can be augmented into an elliptic one.

Nevertheless, it turns out that the zeroth and first cohomology groups of $E_{\theta, A}(M)$ are always finite dimensional. To see this, we first note that the zeroth cohomology $H^0(E_{\theta, A}(M))$ is the Lie algebra of the Stabilizer group $\Stab(A) \subseteq SU(r)$, which is automatically finite dimensional. The first cohomology $H^1(E_{\theta, A}(M))$ is the Zariski tangent space of $\M_{\theta}(M, P)$ at $[A]$ (cf. \cite[Chapter 4]{DK90}), and to see why it is finite dimensional we introduce the following operator 
\begin{equation}
\begin{split}
\delta^{\theta}_A: \Omega^1(M, \g_P) & \longrightarrow \Omega^{n-2}(M, \g_P). \\
a & \longmapsto \theta \wedge d_A a + \star d_Aa 
\end{split}
\end{equation}
Since $\lambda_\theta=-1$, it is clear that for $a \in \Omega^1(M, \g_P)$ we have
\[
\pi_{\theta^{\perp}}(d_A a) = 0 \quad \Longleftrightarrow \quad \delta^{\theta}_A a = 0. 
\]

\begin{lem}
Under the set-up above, $\dim H^1(E_{\theta, A}(M)) < \infty$. Moreover, the complex $E_{\theta,A}(M)$ is elliptic if and only if
\begin{equation*}
    \textup{rank}_\mathbb{R}\;\Lambda^2_{\lambda_\theta}T^\ast M =  \binom{n-1}{2}.
\end{equation*}
\end{lem}

\begin{proof}
We note that $H^1(E_{\theta, A}(M)) = \ker \delta^{\theta}_A / \im d_A$. Appealing to standard elliptic theory \cite{LM89}, it suffices to prove that the symbol sequence of 
\[
\Omega^0(M, \g_P) \xrightarrow{d_A} \Omega^1(M, \g_P) \xrightarrow{\delta^{\theta}_A} \Omega^2(M, \g_P)
\]
is exact at $\Omega^1(M, \g_P)$. We denote the symbol of $\delta^{\theta}_A$ by $\sigma_{\delta}$. Let $\xi \in T^*_pM$ be a non-zero covector on $M$. Then the symbol of $\delta^{\theta}_A$ takes the form 
\[
\sigma_{\delta}(\xi): \eta \longmapsto \theta \wedge (\xi \wedge \eta) + \star (\xi \wedge \eta), \quad \eta \in T^*_pM \otimes \mathfrak{su}(r). 
\]
Suppose $\sigma_{\delta}(\xi)|_{\eta} = 0$. We compute that  
\[
| \xi \wedge \eta |^2 = -2r\tr\left((\xi \wedge \eta) \wedge \star( \xi \wedge \eta) \right)= 2r\tr\left( \theta \wedge (\xi \wedge \eta) \wedge (\xi \wedge \eta) \right)= 0. 
\]
Since $\xi \neq 0$, $\xi \wedge \eta = 0$ implies that $\eta = \xi \otimes \alpha$ for some $\alpha \in \mathfrak{su}(r)$, which means that $\eta$ lies in the image of the symbol $\sigma_{d_A}(\xi)$. The second assertion follows from a dimension count and the fact that the kernel the map $T^\ast_p M\rightarrow \Lambda^2 T^\ast_p M$ given by $\eta\mapsto \xi\wedge \eta$ for $\xi \neq 0$ is precisely the span of $\xi$.
\end{proof}

The failure of ellipticity for $E_{\theta, A}(M)$ tells us that one might not be able to find a rank-finite obstruction sheaf for the moduli space $\M_{\theta}(M, P)$.

\subsection{Anti-Self-Dual Connections on $4$-Manifolds}

Let $X$ be an oriented, Riemannian $4$-manifold. We choose the zero form $\theta \in \Omega^0(X)$ to be the constant $1$ function. Then $\mathcal{Q}_{\theta}$ is simply the Hodge star operator which decomposes the alternating bundle $\Lambda^2T^*X$ into the self-dual and anti-self-dual sub-bundles:
\begin{equation}
\Lambda^2T^*X = \Lambda^+T^*X \oplus \Lambda^-T^*X
\end{equation}
specified by the eigen-decomposition of the Hodge star operator
\begin{equation}
\star \omega = \pm \omega, \quad \omega \in \Lambda^{\pm}T^*X. 
\end{equation}
We write $\omega^{\pm}$ for the projection of $\omega$ onto the subspaces $\Lambda^{\pm}$ respectively. 

\begin{dfn}
Let $\pi_X: P_X \to X$ be a principal $SU(r)$-bundle. An $SU(r)$-connection $B$ on $P_X$ is called anti-self-dual (ASD) if $F_B \in \Omega^-(X, \g_{P_X})$. 
\end{dfn}

We denote the moduli space of ASD connections on $P_X$ by 
\begin{equation}
\mathcal{N}(X, P_X):= \{ B \text{ is an $SU(r)$-connection on $P_X$} : F_B^+ = 0\}/\G_{P_X}. 
\end{equation}
The deformation complex of ASD connections on $4$-manifolds is an elliptic complex:
\begin{equation}\tag{$E_{+, B}(X)$}
0 \to \Omega^0(X, \g_{P_X}) \xrightarrow{-d_B} \Omega^1(X, \g_{P_X}) \xrightarrow{d^+_B} \Omega^+(X, \g_{P_X}) \to 0. 
\end{equation}
The second cohomology $H^2(E_{+, B}(X))$ is commonly referred to as the obstruction space at $[B]$. The local structure near $[B] \in \mathcal{N}(X, P_X)$ is given by the $\Stab(B)$-quotient of the zero set of an obstruction map $\mathfrak{o}_B: H^1(E_{+, B}(X)) \to H^2(E_{+, B}(X))$ as in \cite{DK90}. 

When $X$ is compact, the moduli space $\mathcal{N}(X, P_X)$ admits a well-behaved compactification $\overline{\mathcal{N}}(X, P_X)$ called its Uhlenbeck compactification. Given a positive integer $s \in \N_+$, we write $P_X(s)$ for the $SU(r)$-bundle specified by 
\[
c_2(P_X(s)) = c_2(P_X) - s\PD([X]) \in H^4(X; \Z). 
\] 
An ideal instanton is a pair $([B_o], \pmb{x}) \in \A_{P_X(s)} \times \Sym^s(X)$ with $F^+_{B_o} = 0$. A sequence of ASD instantons $\{[B_n]\}_{n \geq 1}$ is said to converge to an ideal instanton $([B_o], \pmb{x})$ if 
\begin{enumerate}
\item the curvature density $|F_{B_n}|^2$ converges weakly to $|F_{B_o}|^2 + 8\pi^2 \sum_{x \in \pmb{x}} \delta_{x}$, where $\delta_x$ is the delta distribution concentrated at $x$;
\item there exist bundle isomorphisms $v_n: P_X(s)|_{X\backslash \pmb{x}} \to P_X|_{X \backslash \pmb{x}}$ so that $v_n^*B_n$ converges to $B_o$ in $C^{\infty}(K_X)$ for any compact subset $K_X \subseteq X$. 
\end{enumerate} 
The Uhlenbeck's compactness result then states that the moduli space
\begin{equation}
\overline{\mathcal{N}}(X, P_X) := \mathcal{N}(X, P_X)\cup \bigcup_{s=1}^{c_2(P_X)} \mathcal{N}(X, P_X(s)) \times \Sym^s(X)
\end{equation}
of ideal instantons is compact. 

Of particular interest to our paper are the hyperk\"ahler 4-manifolds. Over such manifolds, one can find a triple of K\"ahler forms $\omega_1, \omega_2, \omega_3$ that pointwise span the self-dual bundle $\Lambda^+T^*X$. Thus a connection $B$ is ASD if and only if 
\begin{equation}
F_B \wedge \omega_i = 0, \quad \forall i = 1, 2, 3. 
\end{equation}
If $X$ is compact then it is biholomorphic either to a 4-torus or to a K3 surface. 

\subsection{Hermitian Yang--Mills Connections}

Let $Z$ be a compact K\"ahler manifold of complex dimension $m$. Denote by $\omega$ the K\"ahler form over $Z$. The complex structure on $Z$ equips the complex alternating tensors with a bigrading. To distinguish with the bigrading given by a product structure, we write $\Lambda^{(s,t)} T^*_{\C} Z$ for the bundle of complex $(s, t)$-forms.

\begin{dfn}
Let $\pi_Z: P_Z \to Z$ be a principal $U(r)$-bundle over a K\"ahler manifold $(Z, \omega)$. A unitary connection $B$ on $P_Z$ is called a Hermitian Yang--Mills (HYM) connection if 
\begin{equation}
F^{(2,0)}_B = F^{(0,2)}_B = 0 \qquad iF^{(1,1)}_B \lrcorner \; \omega = \mu \id,
\end{equation}
where $\mu$ is the constant given by 
\begin{equation}
\mu = \frac{2\pi}{(m-1)! \Vol(Z)} \cdot \frac{\deg P_Z}{r},
\end{equation}
and $\deg P_Z = \langle c_1(P_Z) \smile \omega^{m-1}, [Z] \rangle$ is called the degree of $P_Z$. We write $\M_{\omega}(Z, P_Z)$ for the moduli space of Hermitian Yang--Mills connections on $P_Z$. 
\end{dfn}

Let us write $\widehat{F}_B: = iF^{(1,1)}_B \lrcorner \; \omega$ for the section of $\g_P$ obtained from $i F_B$ by contraction with the K\"ahler form $\omega$. The metric on on $\mathfrak{u}(r)$ induces a metric on $\g_{P_Z}$-valued forms via the formula
\[
\langle a, b \rangle : = -2r\tr(ab), \quad \forall a, b \in \mathfrak{u}(r),
\]
which extends the aforementioned metric on $\su(r)$, and decomposes $\mathfrak{u}(r)$ orthogonally as $\mathfrak{u}(r) = i\R \oplus \su(r)$. We note that the curvature $F_B$ satisfies the following pointwise identity \cite[p.405]{IN90}:
\begin{equation}\label{e2.20}
\frac{2r}{(m-2)!} \tr(F_B \wedge F_B) \wedge \omega^{m-2} = \left(|F^{(1,1)}_B|^2 - 2|F_B^{(2,0)}|^2 - \frac{1}{m} |\widehat{F}_B|^2 \right) \vol. 
\end{equation}
Chern--Weil theory tells us that $\tr(F_B \wedge F_B)$ represents the class $8\pi^2 \kappa(P_Z)$, where $\kappa(P_Z) =  c_2(P_Z) - (1/2) c_1^2(P_Z)$. Thus we conclude that the Yang--Mills functional satisfies 
\begin{equation}\label{e2.21}
\begin{split}
\int_Z |F_B|^2 \vol & = \int_Z \frac{2r}{(m-2)!}\tr(F_B \wedge F_B) \wedge \omega^{m-2} + \int_Z 4|F_B^{(2,0)}|^2 +\frac{1}{m} |\widehat{F}_B|^2 \vol \\
& = \frac{16r\pi^2}{(m-2)!} \langle \kappa(P_Z) \smile [\omega]^{m-2}, [Z] \rangle  + \frac{4r\pi^2}{m!} \langle c_1(P_Z) \smile [\omega]^{m-1}, [Z] \rangle\\
&+ \int_Z 4|F_B^{(2,0)}|^2 \vol + \frac{1}{m} \int_Z |\widehat{F}_B - \mu\id|^2 \vol.
\end{split}
\end{equation}
In particular, Hermitian Yang--Mills connections are absolute minimizers of the Yang--Mills functional over a K\"ahler manifold.

When the structure group of $P_Z$ is $SU(r)$, it is straightforward to verify that Hermitian Yang--Mills connections are $\theta$-instantons with $\theta = \omega^{m-2}/(m-2)!$. 

\subsection{$G_2$-Instantons}

Identifying $\R^7$ with the imaginary part $\Ima\mathbb{O}$ of the octonions $\mathbb{O}$, one can define a cross product $u \times v:= \Ima uv$ on the oriented $7$-dimensional Euclidean space. The Lie group $G_2$ is defined as the subgroup of $SO(7)$ that preserves this cross product. Such information can be encoded in a $3$-form $\varphi_0 \in \Lambda^3 (\R^7)^*$ defined by  
\begin{equation}\label{e2.22}
\varphi_0(u, v, w):= \langle u \times v, w \rangle, \quad \forall u, v, w \in \R^7.
\end{equation}
From the formula \eqref{e2.22}, we see that both the inner product and the cross product are determined by the $3$-form. 

Let $M$ be a closed, oriented, smooth $7$-manifold. We write $\mathscr{P}^3(M) \subseteq \Lambda^3T^*M$ for the sub-bundle consisting of alternating tensors that are isomorphic to $\varphi_0$ via orientation-preserving isomorphisms on the tangent spaces of $M$. The fibre of $\mathscr{P}^3_p(M)$ at a point $p \in M$ is then identified with $GL^+(7, \R)/G_2$. We call sections of $\mathscr{P}^3(M)$ positive $3$-forms. In particular, a section $\varphi$ of $\mathscr{P}^3(M)$ induces a metric $g_{\varphi}$ and volume form $\vol_{\varphi}$ on $M$ by the formula
\begin{equation}\label{e2.23}
\iota_u \varphi \wedge \iota_v \varphi \wedge \varphi = 6g_{\varphi}(u, v) \vol_{\varphi}, \quad \forall u, v \in TM. 
\end{equation} 

\begin{dfn}
A torsion-free $G_2$-structure $\varphi$ on an oriented $7$-manifold $M$ is a smooth positive $3$-form $\varphi \in C^{\infty}(M, \mathscr{P}^3(M))$ satisfying 
\[
d\varphi = 0 \qquad d\star_{\varphi} \varphi = 0,
\]
where $\star_{\varphi}$ is the Hodge star given by the metric $g_{\varphi}$ on $M$ induced from $\varphi$. 
\end{dfn}

It is well-known (cf. \cite[Proposition 10.1.3]{J00}) that a section $\varphi$ defines a torsion-free $G_2$-structure if and only if the holonomy group $\Hol(\nabla^{\LC}_{g_{\varphi}})$ of the Levi--Civita connection corresponding to the induced metric $g_{\varphi}$ is contained in $G_2$.

Given a torsion-free $G_2$-manifold $(M, \varphi)$, the corresponding alternating tensor bundle $\Lambda^*T^*M$ decomposes parallelly into fibrewise $G_2$-invariant subspaces. The part closely related to our discussion is the second exterior power
\[
\Lambda^2T^*M = \Lambda^2_7 \oplus \Lambda^2_{14}
\]
given by the fundamental representation of $G_2$ on $\R^7$ and the adjoint representation of $G_2$ on its Lie algebra $\g_2$. Explicitly, we have 
\begin{equation}\label{e2.24}
\Lambda^2_7=\left\{ \omega \in \Lambda^2: \star_{\varphi}(\varphi \wedge \omega)  = 2\omega \right\}  \qquad \Lambda^2_{14} = \left\{ \omega \in \Lambda^2 : \star_{\varphi}(\varphi \wedge \omega) = - \omega \right\}. 
\end{equation}

Note that $\Lambda^2_7$ and $\Lambda^2_{14}$ are respectively the $2$- and $-1$-eigenspaces of the operator $\star_{\varphi}(\varphi \wedge -)$. Thus we obtain the following identity 
\begin{equation}\label{e2.25}
\omega \wedge \omega \wedge \varphi = \left(2|\pi_7(\omega)|^2 - |\pi_{14}(\omega)|^2 \right) \vol_{\varphi}, \quad \forall \: \omega \in \Lambda^2T^*M,
\end{equation}
where $\vol_{\varphi}$ is the volume form induced by the metric $g_{\varphi}$. 

Let $\pi: P \to M$ be a principal $SU(r)$-bundle, and $A_0$ a reference connection on $P$. On the space $\A_P$ of smooth $SU(r)$-connections, we have the $G_2$-version of the Chern--Simons functional:
\begin{equation}
\cs_{\varphi}(A):= -2r \int_M \left(\frac{1}{2} \tr(a \wedge F_A) + \frac{1}{3}\tr(a \wedge a \wedge a) \right) \wedge \star_{\varphi} \varphi,
\end{equation}
where $a = A - A_0$. The gradient of the Chern--Simons functional can be computed to be
\begin{equation}
\cs_{\varphi}(A) = \star_{\varphi}(F_A \wedge \star_{\varphi} \varphi).
\end{equation}
Denote by $\G_P$ the gauge group of $P$. The gradient $\cs_{\varphi}(A)$ is $\G_P$-equivariant. 

\begin{dfn}
Let $P$ be a principal $SU(r)$-bundle over a compact, oriented, smooth, torsion-free $G_2$-manifold $(M, \varphi)$. A connection $A$ on $P$ is called a $G_2$-instanton if 
\begin{equation}\label{e2.28}
F_A \wedge \star_{\varphi} \varphi = 0.
\end{equation}
The moduli space of $G_2$-instantons on $P$ is defined as 
\begin{equation}
\M_{\varphi}(M, P):= \left\{A \in \A_P : F_A \wedge \star_{\varphi} \varphi = 0 \right\}/\G_P. 
\end{equation}
\end{dfn}

After unraveling the defining equation \eqref{e2.25}, we have an equivalent characterization of $G_2$ instantons (see for instance \cite{S89}):
\begin{equation}\label{e2.30}
F_A \wedge \star_{\varphi} \varphi = 0 \quad \Longleftrightarrow \quad \varphi \wedge F_A = -\star_{\varphi} F_A.
\end{equation}
Thus $G_2$-instantons are $\theta$-instantons by taking $\theta = \varphi$. From \eqref{e2.31}, we see that the eigenvalues of $\mathcal{Q}_{\theta}$ consist of $-1$ and $2$. \autoref{p2.3} then tells us that $G_2$-instantons minimize the Yang--Mills functional. 

Appealing to \eqref{e2.25}, the deformation complex of $G_2$-instantons can be augmented to an elliptic complex:
\begin{equation}\tag{$E_{\varphi, A}(M)$}
0 \to \Omega^0 \xrightarrow{-d_A} \Omega^1\xrightarrow{\star_{\varphi} \varphi \wedge d_A} \Omega^6 \xrightarrow{d_A} \Omega^7 \to 0
\end{equation}

\subsection{$\Spin(7)$-Instantons}

After identifying the oriented $8$-dimensional Euclidean space $\R^8$ with the octonions $\mathbb{O}$, we get a triple cross product defined by 
\[
u \times v \times w := \frac{1}{2} \left((u\bar{v})w - (w\bar{v})u \right),
\]
from which we get a self-dual $4$-form $\Phi_0 \in \Lambda^4(\R^8)^*$ via the formula
\begin{equation}\label{e2.31}
\Phi_0(u, v, w, x):= \langle u \times v \times w, x\rangle, \quad \forall u, v, w, x \in \R^8.
\end{equation}
The compact Lie group $\Spin(7)$ is the subgroup of $SO(8)$ that preserves the $4$-form $\Phi_0$. 

Let $M$ be a closed, oriented, smooth $8$-manifold. We denote by $\mathscr{S}^4(M) \subseteq \Lambda^4T^*M$ the sub-bundle consisting of alternating tensors that are isomorphic to $\Phi_0$ via orientation-preserving isomorphisms on the tangent spaces of $M$. The fibre $\mathscr{S}^4_p(M)$ at $p \in M$ is identified with $GL^+(8, \R)/\Spin(7)$. Sections of $\mathscr{S}^4(M)$ are called admissible $4$-forms. A section $\Phi$ of $\mathscr{S}^4(M)$ also induces a metric $g_{\Phi}$ and volume form $\vol_{\Phi}$ on $M$ via the formula
\begin{equation}\label{e2.32}
\iota_v\iota_u\Phi \wedge \iota_w\iota_u\Phi \wedge \Phi = 6\left(|u|^2_{g_{\Phi}}\langle v, w \rangle_{g_{\Phi}} - \langle v, u \rangle_{g_{\Phi}}\langle u, w\rangle_{g_{\Phi}} \right)\vol_{\Phi}.
\end{equation}

\begin{dfn}
A torsion-free $\Spin(7)$-structure on an oriented $8$-manifold $M$ is a smooth, admissible $4$-form $\Phi \in C^{\infty}(M, \mathscr{S}^4(M))$ satisfying $d\Phi = 0$.
\end{dfn}

By \cite[Proposition 10.5.3]{J00}, an admissible $4$-form $\Phi$ is closed if and only if $\Phi$ is parallel with respect to the Levi--Civita connection of the induced metric $g_{\Phi}$, which reduces the holonomy group to $\Hol(g_{\Phi}) \subseteq \Spin(7)$. 

A torsion-free $\Spin(7)$-structure also decomposes the alternating tensor bundle $\Lambda^*T^*M$ into $\Spin(7)$-invariant sub-bundles. Of our particular interest is the decomposition of the $2$-forms
\[
\Lambda^2T^*M = \Lambda^2_7 \oplus \Lambda^2_{21},
\]
where $\Lambda^2_7$ and $\Lambda^2_{21}$ are sub-bundles of rank $7$ and $21$ respectively, and given explicitly as 
\[
\Lambda^2_7 = \{ \omega \in \Lambda^2: \star_{\Phi}(\Phi \wedge \omega) = 3\omega \} \qquad \Lambda^2_{21} = \{ \omega \in \Lambda^2: \star_{\Phi}(\Phi \wedge \omega) = -\omega \}. 
\]

\begin{dfn}
Let $(M, \Phi)$ be a torsion-free $\Spin(7)$-manifold and $\pi: P \to M$ a principal $SU(r)$-bundle. A connection $A$ on $P$ is called a $\Spin(7)$-instanton if 
\[
\star_{\Phi}( F_A \wedge \Phi)  = -F_A.
\]
\end{dfn}

From the definition, we read that $\Spin(7)$-instantons are $\theta$-instantons with $\theta = \Phi$. The eigenvalues of $\mathcal{Q}_{\theta}$ are $-1$ and $3$. Thus $\Spin(7)$-instantons minimize the Yang--Mills functional as well. The deformation complex of $\Spin(7)$-instantons is elliptic (cf. \cite{DS11}):
\begin{equation}\tag{$E_{\Phi, A}(M)$}
0 \to \Omega^0 \xrightarrow{-d_A} \Omega^1 \xrightarrow{\pi_7 \circ d_A} \Omega^2_7 \to 0,
\end{equation}
where $\pi_7: \Lambda^2 \to \Lambda^2_7$ is the orthogonal projection.

\section{Curvature and Parallelism}\label{s3}

Let $G$ be a compact Lie subgroup of some general linear group. We consider a principal $G$-bundle $\pi: P\rightarrow Y\times X$ over a product manifold. The product structure then induces a bigrading on differential forms defined over $Y \times X$. Given a $G$-connection $A$ on $P$, we shall analyze the system of equations
\begin{equation}
    F^{2,0}_A=0 \qquad \textup{and} \qquad  F^{1,1}_A=0.
\end{equation}
Regarding $P$ as a family of $G$-bundles over $X$, these equations correspond to the parallel condition for the restricted connections $A|_{\{y\} \times X}$. We shall establish stability and integrability results for such family, and prove \autoref{t1.6}. 

\subsection{Bundles over Product Manifolds}\label{s3.1}

Let $Y$ and $X$ be two closed, smooth, positive-dimensional Riemannian manifolds. We write $M = Y \times X$ for the product manifold and $\pi_1: M \to Y$, $\pi_2: M \to X$ for the factor projections. We regard $M$ as a fibre bundle over the base $Y$. For any $y \in Y$, we write $X_y:= \{y\} \times X$ for the fibre over $y$. 

Recall that $G$ is a compact linear group. We denote by $\g:= \Lie(G)$ the Lie algebra of $G$. Let $\pi: P \to M$ be a principal $G$-bundle. We write $A$ for a smooth $G$-connection on $P$ and $\g_P := P \times_{\ad} \g$ for the adjoint bundle. The space $\A_P$ of smooth $G$-connections forms an affine space modeled on $\Omega^1(M, \g_P)$. 

Since $M = Y \times X$ is a product manifold, the space of alternating tensors inherits a bigrading 
\[
\Lambda^l T^*M = \bigoplus_{s + t = l} \Lambda^s \pi^*_1 T^*Y \otimes \Lambda^t \pi^*_2 T^*X. 
\]
The $\g_P$-valued forms on $M$ decompose accordingly as 
\[
\Omega^l(M, \g_P) := \bigoplus_{s+t=l} \Omega^{s, t}(M, \g_P). 
\]
Given a $G$-connection $A$, the coupled exterior derivative decomposes as 
\[
d_A = d_{H, A} + d_{f, A}: \Omega^{s, t} \longrightarrow \Omega^{s+1, t} \oplus \Omega^{s, t+1}.  
\]
with respect to this bigrading and the curvature $F_A$ can be written as
\[
F_A = F^{2,0}_A + F^{1,1}_A + F^{0,2}_A. 
\]

We write $P_y = P|_{X_y}$ for the restriction of the bundle $P$ to the fibre $X_y$ and $B_y := A|_{P_y}$ for the restriction of the connection $A$. We fix a reference connection $A_0 \in \A_P$. The difference between $A_0$ and $A$ can be written as 
\[
A-A_0 = a^{1,0} + a^{0,1} \in \Omega^{1,0}(M, \g_P) \oplus \Omega^{0,1}(M, \g_P). 
\]
We write  $a^{0,1}|_{P_y} = b_y \in \Omega^1(X, P_y)$ for the restriction of the $(0,1)$-component. Let $(U, y_1, ..., y_m)$ be a local chart for $Y$. Over $U \times X$, the $(1,0)$-component $a^{1,0}$ can be written as 
\[
a^{1,0}:= \sum_{i=1}^m \sigma_idy_i, \quad \text{ where } \sigma_i(y, -) \in \Omega^0(X, P_y). 
\]
Then each component of the curvature $F_A$ takes the following form. Writing 
\begin{equation}
F^{ij}_{A_0, \sigma} = \nabla_{A_0, \partial_{y_i}} \sigma_j - \nabla_{A_0, \partial_{y_j}} \sigma_i + [\sigma_i, \sigma_j],
\end{equation}
we have 
\begin{equation}\label{e3.3}
F_A^{2,0} - F_{A_0}^{2,0} = d_{H, A_0}a^{1,0} + a^{1,0} \wedge a^{1,0}  = \sum_{i<j} F^{ij}_{A_0, \sigma} dy_i \wedge dy_j. 
\end{equation}
As the for $(1,1)$-component, we have 
\begin{equation}\label{e3.4}
F^{1,1}_A - F^{1,1}_{A_0} = \sum_{i=1}^m dy_i \wedge (\nabla_{A_0, \partial_{y_i}} b_y - d_{B_y} \sigma_i). 
\end{equation}
Finally, the $(0, 2)$-component is 
\begin{equation}
F^{0,2}_A = F_{B_y}. 
\end{equation}

Since the bundle $P$ is defined over a fibre bundle $M = Y \times X$, we can shift our viewpoint to consider the associated family bundle $\underline{\pi}: \underline{P} \to Y$ whose fibres consist of $P_y$. Note that the total space $\underline{P}$ is the same as that of the original bundle $P$, and the bundle projection map is just the composition $\underline{\pi} = \pi_1 \circ \pi$. We fix a bundle $P_X \to X$ that is isomorphic to the fibres $P_y$. Then the structure group of $\underline{P}$ is the gauge group $\G_{P_X}$ of $P_X$. 

A connection $A$ on $P$ induces an Ehresmann connection $\underline{A}$ on $\underline{P}$ as follows. The connection $A$ induces a $G$-invariant decomposition $TP = H_A \oplus VP$ with vertical bundle $VP = \ker (\pi_*:TP \to \pi^*TM)$. Appealing to the isomorphism
\begin{equation*}
    \pi_*: H_A \to \pi^*TM = \underline{\pi}^*TY \oplus (\pi_2 \circ\pi)^*TX,
\end{equation*}
we define 
\begin{equation}\label{e3.6}
H_{\underline{A}}:= \ker (p_2 \circ \pi_*: H_A \to (\pi_2 \circ\pi)^*TX) \subseteq T\underline{P}
\end{equation}
as the horizontal distribution defining the connection $\underline{A}$, where $p_2: \underline{\pi}^*TY \oplus (\pi_2 \circ\pi)^*TX \to (\pi_2 \circ\pi)^*TX$ is projection to the second factor.

\begin{lem}\label{l3.1}
The curvature $F_{\underline{A}}$ of the induced connection is identified canonically with the $(2,0)$-component $F^{2,0}_A$ of the original connection on $P$.
\end{lem}
\begin{proof}
We have a decomposition $T\underline{P} = H_{\underline{A}} \oplus V\underline{P}$ with $V\underline{P} = (\pi_2 \circ \pi)^*TX \oplus VP$. Denote by $\underline{\eta}:TY \to C^{\infty}(\underline{P}, H_{\underline{A}})$ the $H_{\underline{A}}$-horizontal lift in $T\underline{P}$. Given two tangent vectors $v, w \in T_yY$, the curvature is defined by 
\begin{equation}
F_{\underline{A}}|_y(v, w):=[\underline{\eta}(v), \underline{\eta}(w)] - \underline{\eta}([v, w]) \in V\underline{P}|_{P_y},
\end{equation}
where we have used the fact that $\underline{\pi}_* \circ \underline{\eta} = \id$. We note that $(\pi_2 \circ \pi)^*TX|_{P_y} = \pi^*_yTX_y$ and $VP|_{P_y} = \g_{P_y}$, where $\pi_y: P_y \to X_y$ is the projection map. Thus if we further fix $x \in X_y$, the curvature element $F_{\underline{A}}|_{(y, x)}(v, w)$ defines a tangent vector $\alpha \in T_xX_y$ and an element $\xi \in \g_{P_y}$. Then it follows tautologically that  
\begin{equation}
F^{2,0}_A|_{y, x}((v, \alpha), (w, \alpha)) =(j_y)_*\left( [\eta(v, \alpha), \eta(w, \alpha)] - \eta[(v, \alpha), (w, \alpha)] \right)= \xi,
\end{equation}
where $j_y: P_y \times \su(r) \to \g_{P_y}$ is the quotient map, and $\eta: T(Y \times X) \to C^{\infty}(P, H_A)$ is the $H_A$-horizontal lift in $TP$. This identifies $F_{\underline{A}}$ with $F^{2,0}_A$.
\end{proof}

If one wishes to work with principal bundles rather than the fibre bundle $\underline{P}$, one can consider the bundle $\mathbb{P}$ of fibrewise isomorphisms of $\underline{P}$. Then $\mathbb{P}$ is a principal $\G_{P_X}$-bundle. Since the gauge group $\G_{P_X}$ acts effectively on the fibres of $\underline{P}$, the notions of connection and curvature are corresponded properly. One can certainly discuss characteristic classes of such bundles, which come from cohomology of the classifying space $B\G_{P_X}$. Since the cohomology of $B\G_{P_X}$ depends on the manifold $X$, we shall not pursue a general theory here. 

We end this subsection with the following definition.

\begin{dfn}\label{d2.1}
Let $\underline{A}$ be the induced connection on the family bundle $\underline{P}$. We call its holonomy the horizontal holonomy of $A$, denoted by $\Hol^h_A$, which is a subgroup of $\G_{P_X}$ defined up to conjugation. 
\end{dfn}

Once we fix an identification between $P_{y_0}$ and $P_X$, the horizontal holonomy $\Hol^h_A$ becomes a well-defined subgroup of $\G_{P_X}$.

\subsection{Stabilization}\label{s3.2}

Let us choose a basepoint $y_0 \in Y$. We identify the restriction of $P$ to the slice $\{y_0\} \times X$ with a reference bundle $\pi_X: P_X \to X$ whose gauge group is denoted by $\G_{P_X}$. Denote by $L(Y,y_0)$ the space of smooth loops in $Y$ based at $y_0$. Given $\gamma\in L(Y,y_0)$ and a point $p\in P_X$, we get a loop 
\begin{equation}
\begin{split}
\gamma_{\pi(p)}: [0, 1] &\longrightarrow Y \times X \\
t & \longmapsto (\gamma(t), \pi_X(p))
\end{split}
\end{equation}
Let $A$ be a $G$-connection on $P$ whose corresponding horizontal distribution is denoted by $H_A \subseteq TP$. We write $\tilde{\gamma}^A_{\pi(p)}$ the unique $H_A$-horizontal lift of the loop $\gamma_{\pi(p)}$ to $P$ that is based at $p$. 

\begin{dfn}
We define the \textit{horizontal holonomy representation} of the connection $A$ to be the map
\begin{equation}
\begin{split}
    \mathfrak{H}_A: L(Y,y_0) & \longrightarrow \G_{P_X} \\
 \gamma & \longmapsto \left(p \mapsto \tilde{\gamma}^A_{\pi(p)}(1)\right).
\end{split}
\end{equation}
\end{dfn}

To see why $\mathfrak{H}_A(\gamma)$ defines an element in the gauge group $\G_{P_X}$, we first note that $\mathfrak{H}_A(\gamma)$ is an automorphism of the bundle $P_X \to X$ that covers the identity map on $X$. The equivariance of $H_A$ implies that $\tilde{\gamma}^A_{\pi(p\cdot g)} = \tilde{\gamma}^A_{\pi(p)} \cdot g$ for $\forall g \in G$. Thus $\mathfrak{H}_A(\gamma)$ is $G$-equivariant and defines a gauge transformation on $P_X$. 

Recall that the holonomy $\Hol_{\underline{A}}$ of the induced connection on $\underline{P}$ is referred to as the horizontal holonomy of $A$ in \autoref{d2.1}. Our choice of terminology is justified by the following lemma.

\begin{lem}
For any $\gamma \in L(Y, y_0)$, we have 
\[
\Hol_{\underline{A}}(\gamma) = \mathfrak{H}_A(\gamma) \in \G_{P_X}.
\]
In particular, when $F^{2,0}_A = 0$, the horizontal holonomy defines a group homomorphism $\pi_1(Y,y_0)\rightarrow\mathcal{G}_X$.
\end{lem}

\begin{proof}
Given a loop $\gamma \in L(Y, y_0)$, we write $\tilde{\gamma}^{\underline{A}}_p$ for the $\underline{A}$-horizontal lift of $\gamma$ based at $p \in P_X$. We note that $\gamma_{\pi(p)}$ is constant along the $X$-direction. Thus 
\[
p_2 \circ \pi_*\left(\partial_t\tilde{\gamma}^A_{\pi(p)}\right) = 0,
\]
where $p_2 \circ \pi_*$ is the map considered in \eqref{e3.6}. From the construction of $H_{\underline{A}}$, it follows that $\partial_t\tilde{\gamma}^A_{\pi(p)} \in H_{\underline{A}}$. The uniqueness of the horizontal lift then tells us that $\tilde{\gamma}^{\underline{A}}_p = \tilde{\gamma}^A_{\pi(p)}$, which proves the first statement. 

As for the second statement, we note that $F^{2,0}_A= 0 $ implies that $F_{\underline{A}} = 0$ by \autoref{l3.1}. So $\underline{A}$ becomes a flat connection, which gives rise to the claimed homomorphism. 
\end{proof}

Now we proceed to relate the vanishing of the $(1,1)$-component of $F_A$ to the horizontal holonomy. For computational reasons, we work on the associated vector bundle $E=P \times_G \mathbb{V}$, where $\mathbb{V}$ is a vector space on which $G$ acts effectively. Let $J_t =[0, 1]$ and $I_s = (-1, 1)$ be two intervals parametrized by $t$ and $s$, respectively. Given two smooth paths $\alpha: J_t \rightarrow X$ and $\gamma: I_s \rightarrow Y$, we consider the following smooth immersion
\begin{equation}
\begin{split}
    \gamma_\alpha:  J_t \times I_s & \longrightarrow Y\times X \\
(t, s) & \longmapsto (\gamma(t), \alpha(s)). 
\end{split}
\end{equation}
We shall write $\gamma_{\alpha(s)}$ for the path given by $\gamma_{\alpha}(-, s)$ and we assume that $\gamma_{\alpha(s)}(0) = (y_0, \alpha(s))$. Denote by $\Pi^A_{t, s}: E|_{\gamma_{\alpha}(t, s)} \to E|_{\gamma_{\alpha}(1, s)}$ the parallel transport map along the path $\gamma_{\alpha(s)}$ from $\gamma_{\alpha}(t, s)$ to $\gamma_{\alpha}(1, s)$. 

We wish to compute the differential of $\Pi_{t, s}$ with respect to $s$. To formulate this, we momentarily consider the bundle $\mathcal{H} \to I_s$ whose fibre at $s$ is $\Hom(E|_{\gamma_{\alpha}(0, s)}, E|_{\gamma_{\alpha}(1,s)})$. Then the parallel transport map $\Pi^A_{0, s}$ is a section of $\mathcal{H}$. The intrinsic derivative of $\Pi_{0, s}$ is defined as the composition of the differential of $\Pi^A_{0, s}$ with the projection onto the tangential part of the fibres in $T\mathcal{H}$, which we denote by $D_s \Pi^A_{0, s}$. 

\begin{lem}\label{l3.5}
The intrinsic derivative of $\Pi_{0,s}$ with respect to $s$ is given by 
\begin{equation}\label{e3.12}
    D_s|_{s=s_0} \Pi^A_{0, s} = \left ( \int^1_0 \Pi^A_{t, s_0}\circ F_A(\partial_t\gamma_\alpha,\partial_s\gamma_\alpha)\circ (\Pi^A_{t, s_0})^{-1} \; dt\right )\cdot \Pi^A_{0, s_0}, \quad \forall s_0 \in I_s.
\end{equation}
\end{lem}
\begin{proof}
We may work on the pull-back bundle $\gamma_{\alpha}^*E \to J_t \times I_s$ and fix a trivialization of this bundle so that the fibres $E|_{\gamma_{\alpha}(t, s)}$ are identified. Such choices will not affect the computation since the derivative is intrinsically defined. We shall omit writing out the pull-back map $\gamma^*_{\alpha}$ in the following argument when the context is clear. 

Given a vector $\sigma(0, s_0) \in E|_{\gamma_{\alpha}(0, s_0)}$, contractibility of $J_t \times I_s$ enables us to extend the vector as a section $\sigma$ of $\gamma_{\alpha}^*E$ satisfying 
\begin{equation}
    \nabla_{A, \partial_t}\sigma(t, s) \equiv 0, \;\;\; \nabla_{A, \partial_s}\sigma (0, s)\equiv 0
\end{equation}
where $\nabla_A$ is the induced covariant derivative on $J_t \times I_s$. The construction tells us that 
\begin{equation}\label{e3.14}
\Pi^A_{0, s_0} \sigma(0, s_0) = \sigma(1, s_0) \quad \text{ and } \quad (D_s|_{s=s_0} \Pi^A_{0, s}) \sigma(0, s_0) = \nabla_{A, \partial_s}\sigma(1, s_0). 
\end{equation}
From the fact that $\sigma$ is parallel in $t$, it follows that
\begin{equation}
    \nabla_{A, \partial_t}\nabla_{A, \partial_s}\sigma=F_A(\partial_t\gamma_\alpha,\partial_s\gamma_\alpha)\sigma.
\end{equation}
Using this, we compute that
\begin{align*}
\frac{\partial}{\partial t}\Pi^A_{t, s}\nabla_{A, \partial_s}\sigma(t, s) &= \Pi^A_{t, s}\nabla_{A, \partial_t}\nabla_{A, \partial_s}\sigma(t, s)\\
&=\Pi^A_{t, s}F_A(\partial_t\gamma_\alpha,\partial_s\gamma_\alpha)(\Pi^A_{t, s})^{-1}\sigma(1, s)
\end{align*}
On the other hand, since $\sigma(0, s)$ is parallel in $s$ while $\Pi_{1, s}$ is the identity, we can conveniently write
\begin{equation}
\begin{split}
    \nabla_{A, \partial_s}\sigma(1, s_0) & =\Pi^A_{1, s_0}\nabla_{A, \partial_s}\sigma(1, s_0)-\Pi^A_{0, s_0}\nabla_{A, \partial_s}\sigma(0, s_0) \\
    & = \int_0^1 \Pi^A_{t, s_0}F_A(\partial_t\gamma_\alpha,\partial_s\gamma_\alpha)(\Pi^A_{t, s_0})^{-1}dt \cdot \sigma(1, s_0),
\end{split}
\end{equation}
which proves the claimed formula after appealing to \eqref{e3.14}. 
\end{proof}

The computation of the intrinsic derivative helps us rewrite the vanishing condition $F^{1,1}_A = 0$ as a condition on the parallel transport map.

\begin{cor}\label{c3.6}
Let $A$ be a $G$-connection on $E$. Then $F^{1,1}_A = 0$ if and only if the parallel transport map $\Pi^A_{0, s}: E|_{\gamma_{\alpha}(0, s)} \to E|_{\gamma_{\alpha}(1, s)}$ is constant in $s$ for all smooth paths $\gamma: I_s \to Y$, $\alpha: J_t \to X$. 
\end{cor}

\begin{proof}
Suppose $F^{1,1}_A = 0$. We note that $\partial_t \gamma_{\alpha}$ is a section of $\gamma^*_{\alpha}\pi^*_1TY$ and $\partial_s \gamma_{\alpha}$ is a section of $\gamma^*_{\alpha}\pi_2^*TX$. Thus only the $(1,1)$-component of $F_A$ could possibly contribute non-zero terms in the integrand of \eqref{e3.12}. Thus $F^{1,1}_A = 0$ implies that the $D_s\Pi^A_{0,s} = 0$, which means $\Pi^A_{0,s}$ is constant in $s$. 

Conversely suppose $\Pi^A_{0, s}$ is constant for all choices of $\alpha$ and $\gamma$. Given two tangent vectors $v \in T_y Y$ and $w \in T_xX$, we may choose $\gamma$ and $\alpha$ accordingly so that 
\[
\gamma(0) = y, \quad \dot{\gamma}(0) = v, \quad \alpha(0) = x, \quad \dot{\alpha}(0) = w. 
\]
Then \autoref{l3.5} tells us that 
\begin{equation}\label{e3.17}
 \int^1_0 \Pi^A_{t, 0}\circ F_A(\partial_t\gamma_\alpha,\partial_s\gamma_\alpha)\circ (\Pi^A_{t, 0})^{-1} \; dt = 0. 
\end{equation}
For each $t_0 \in (0, 1]$, we can modify $\gamma$ so that $\dot{\gamma}(t) = 0$ for $\forall t \geq t_0$. Then the upper limit of the integral in \eqref{e3.17} becomes $t_0$. Differentiating the integral at $t_0 = 0$ gives us 
\[
\Pi^A_{0, 0} \circ F_A(v, w) \circ (\Pi^A_{0,0})^{-1} = F^{1,1}_A(v, w)= 0.
\] 
Since $v$ and $w$ are chosen arbitrarily, we conclude that $F^{1,1}_A = 0$. 
\end{proof}

With this interpretation in hand, we prove an important property for the horizontal holonomy. 

\begin{prop}\label{p3.7}
Let $A$ be a $G$-connection on $P$ satisfying $F^{1,1}_A = 0$. Then the horizontal holonomy $\mathfrak{H}_A(\gamma)$ stabilizes $A|_{P_{y_0}}$ for every $\gamma \in L(Y, y_0)$. 
\end{prop}

\begin{proof}
Let us write $B = A|_{P_{y_0}}$ which we regard as a connection on $P_X$. Denote by $H_B \subseteq TP_X$ the horizontal distribution induced by the connection $B$. Given $\gamma\in L(Y,y_0)$, we get the differential $D\mathfrak{H}_A(\gamma): TP_X \to TP_X$. Then we know 
\begin{equation}
\mathfrak{H}_A(\gamma) \text{ stabilizes } B \Longleftrightarrow D\mathfrak{H}_A(\gamma)(H_B) \subseteq H_B. 
\end{equation}
Note that $\mathfrak{H}_A(\gamma)$ sends $p \mapsto p\cdot \Hol_A(\gamma_{\pi(p)})$ for $\forall p \in P_X$. So we can view $\mathfrak{H}_A(\gamma) : P_X\rightarrow P_X$ as the composition of the two maps
\begin{equation}
    \mu : P_X\times G \rightarrow P_X, \quad \mu(p,g)=p\cdot g
\end{equation}
\begin{equation}
    \nu: P_X\rightarrow P_X\times G, \quad \nu(p)=\left(p, \Hol_A(\gamma_{\pi(p)})\right)
\end{equation}
so that $D\mathfrak{H}_A(\gamma)  = D\mu \circ D\nu$. We write $\mathfrak{r}_g: P_X \rightarrow P_X$ and $\mathfrak{l}_g: G \rightarrow G$ are right and left multiplication by $g\in G$, respectively. Then the differential $D\mu$ is given by 
\begin{equation}
    D\mu|_{(p,g)}(u, \xi)=D\mathfrak{r}_g|_p(u) + \xi^*, \quad u \in T_pP_X, \; \xi \in \g,
\end{equation}
where $\xi^\ast\in T_{p\cdot g}P_X$ is the tangent vector determined by $\xi$ via the formula
\begin{equation}
    \xi^\ast = \frac{d}{dt}\bigg\rvert_{t=0}\left(pg\cdot e^{t D\mathfrak{l}_{g^{-1}}(\xi)}\right).
\end{equation}
For the convenience of notation, we momentarily write $\tilde{\nu}(p) = \Hol_A(\gamma_{\pi(p)})$ so that $\nu = (\id, \tilde{\nu})$. Then the chain rule enables us to compute $D\mathfrak{H}_A(\gamma)$ as 
\begin{equation}
D\mathfrak{H}_A(\gamma)|_p(u) = D\mu|_{(p, \tilde{\nu}(p))}(u, D\tilde{\nu}|_p u) = D\mathfrak{r}_{\tilde{\nu}(p)}|_p u + (D\tilde{\nu}|_p u)^*, \quad \forall u \in T_pP_X. 
\end{equation} 

With respect to the splitting $TP_X=H_B\oplus VP_X$, denote by $\pi_V: TP_X\rightarrow VP_X$ the vertical projection. From the $G$-equivariance of $H_B$, it follows that
\begin{equation}\label{e3.24}
\pi_V\circ D\mathfrak{H}_A(\gamma)|_p (u) =(D\tilde{\nu}|_p u)^*.
\end{equation}
To compute $D\tilde{\nu}|_p u$, we choose a curve $\alpha: (-\epsilon, \epsilon) \to X$ satisfying 
\[
\alpha(0) = p \qquad \dot{\alpha}(0) = D\pi_X(u) \in T_{\pi(p)}X
\] 
for some small $\epsilon > 0$. Denote by $\tilde{\alpha}$ the $H_B$-horizontal lift in $P_X$ based at $p$. Then \autoref{c3.6} implies that 
\begin{equation}
    D\tilde{\nu}|_pu = \frac{d}{dt}\bigg\rvert_{s=0}\tilde{\nu}(\tilde{\alpha}(s)) = \frac{d}{ds}\bigg\rvert_{s=0}\Hol_A(\gamma_{\pi(\tilde{\alpha}(s))})=0\in V_pP_X. 
\end{equation}
Combining with \eqref{e3.17}, we conclude that $D\mathfrak{H}_A(\gamma)|_p u \in H_B$ for all $u \in H_B$, as desired.
\end{proof}

\subsection{Integrability}

In this subsection, we switch to work with connection $1$-forms rather than horizontal distributions. We wish to argue that the vanishing of the $(2,0)$-component of the curvature is an integrability condition. Then we can identify connections once they have the same holonomy. 

Let $A_0$ be a reference connection on a principal $G$-bundle $P \to Y \times X$ and suppose that
\begin{equation}
F^{2,0}_{A_0} = F^{1,1}_{A_0} = 0. 
\end{equation}
Recall that locally over $U \times X$ with respect to some chart $(U, y_i) \subseteq Y$, we can write any other connection $A$ on $P$ as 
\begin{equation}
A = A_0 + b_y + \sum_{i=1}^m \sigma_i dy_i, \quad \text{ where } b_y \in \Omega^1(X, P_y), \; \sigma_i(y, -) \in \Omega^0(X, P_y). 
\end{equation}
Then from \eqref{e3.3} we read that 
\begin{equation}\label{e3.28}
F^{2,0}_A = 0 \quad \Longleftrightarrow \quad \nabla_{A_0, \partial_{y_i}}\sigma_j - \nabla_{A_0, \partial_{y_j}} \sigma_i + [\sigma_i, \sigma_j] = 0,
\end{equation}
and from \eqref{e3.4} we read that 
\begin{equation}\label{e3.29}
F^{1,1}_A = 0 \quad \Longleftrightarrow \quad \nabla_{A_0, \partial_{y_i}} b_y - d_{f, A_0}\sigma_i - [b_y, \sigma_i]  = 0. 
\end{equation}
The following result realizes \eqref{e3.28} as an integrability condition. 

\begin{prop}\label{p3.8}
Let $P \to Y \times X$ be a principal $G$-bundle. Suppose $A_0$ and $A$ are two smooth connections on $P$ satisfying 
\[
F^{2,0}_{A_0} = F^{2,0}_A = 0 \qquad F^{1,1}_{A_0} = F^{1,1}_A = 0.
\]
Then for any simply-connected domain $\Omega \subseteq Y$ containing the basepoint $y_0$, there exists a unique smooth gauge transformation $u$ on $P|_{\Omega \times X}$ satisfying 
\begin{equation*}
(1) \; u|_{P_{y_0}} = \id \qquad (2)\; u \cdot A - A_0 = \tilde{b}_y \qquad (3)\; d_{H, A_0} \tilde{b}_y = 0 
\end{equation*}
for some smooth family $\tilde{b}_y \in \Omega^1(X, \g_{P_y})$ with $y \in \Omega$. 
\end{prop}

\begin{proof}
We first prove the uniqueness. Suppose we have two gauge transformations $u_1$ and $u_2$ satisfying the desired properties. Let us write $A_i = u_i \cdot A$, $i=1, 2$, and $u' = u_2u^{-1}_1$. Then $u' \cdot A_1 = A_2$. Since both the difference $A_1 - A_0$ and $A_2 - A_0$ have vanishing $(1,0)$-components, we see that $d_{H, A_0} u' = 0$, i.e. $u'$ is $A_0$-parallel along the $Y$-direction. Since $F^{2,0}_{A_0} = 0$ and $\Omega$ is simply-connected, $u'$ is determined (via parallel transport) by its restriction on a single slice $P_y$. Due to the requirement that $u|_{P_{y_0}} = \id$, we see that $u' = \id$. Thus $u_1 = u_2$. 

Now we establish the existence. Given $y \in \Omega$, we construct $u_y \in \Aut(P_y)$ as follows. Let us pick any path $\gamma$ in $\Omega$ connecting $y_0$ to $y$. We identify a tubular neighborhood $\nu(\gamma)$ of $\gamma$ with $[0, 1] \times I^{m-1}_{\epsilon}$ so that $y_0 = (0, 0, ..., 0)$ and $y=(1, 0, ...., 0)$, where $I_{\epsilon} = (-\epsilon, \epsilon)$ is the open interval of length $2\epsilon > 0$. We denote the coordinates on $\nu(\gamma)$ by $(y_1, ..., y_m)$. Then we can write
\[
A = A_0 + b_y + \sum_{i=1}^m \sigma_i dy_i
\]
as before. We shall construct the gauge transformation $u$ on $\nu(\gamma) \times X$ by inductively annihilating the connection forms $\sigma_i$.

The action of $u$ on $A$ can be computed as 
\begin{equation}\label{e3.30}
 u\cdot A = A_0 + b_y + \sum_i \sigma_idy_i - \left(d_{A_0}u + [b_y, u] + \sum_i \left(\nabla_{A_0, \partial_{y_i}} u + [\sigma_i, u] \right) dy_i \right)u^{-1}. 
\end{equation}
Thus it suffices to solve the system of differential equations 
\begin{equation}
\nabla_{A_0, \partial_{y_i}}u = u\sigma_i, \quad i = 1, ..., m. 
\end{equation}
We solve the system inductively on $i$. Starting with $i=1$, we solve the ODE 
\[
\nabla_{A_0, \partial_{y_1}} u = u\sigma_1 \qquad u(0, y_2, ..., y_m) = \id
\]
with respect to fixed $y_2, ..., y_m$ to get a unique smooth gauge transformation $u_1$. Let us write 
\[
u_1 \cdot A = A_0 + b_y^{(1)} + \sum_i \sigma_i^{(1)} dy_i. 
\]
Then \eqref{e3.30} tells us that $\sigma_1^{(1)} = 0$. Suppose inductively we have arrived at the $j$-th step where a gauge transformation $u_j$ has been constructed so that 
\[
u_j \cdot A = A_0 + b^{(j)}_y + \sum_{i=j+1}^m \sigma_i^{(j)}dy_i.
\]
Since $u_j \cdot A$ satisfies the condition $F^{2,0}_{u_j\cdot A}=0$, \eqref{e3.28} tells us 
\begin{equation}\label{e3.32}
\nabla_{A_0, \partial_{y_i}} \sigma^{(j)}_k = 0 \quad \forall i \leq j.  
\end{equation}
Now we solve the ODE
\[
\nabla_{A_0, \partial_{y_{j+1}}} u = u\sigma^{(j)}_{j+1} \qquad u(y_1, ...,y_j, 0, y_{j+2}, ..., y_m) = \id
\]
to get a unique smooth gauge transformation $u'_{j+1}$. Let $u_{j+1} = u'_{j+1} u_j$, and write 
\[
u_{j+1} \cdot A = A_0 + b^{(j+1)}_y + \sum_i \sigma_i^{(j+1)}dy_i.
\]
Then since $[\partial_{y_i}, \partial_{y_j}] = 0$, \eqref{e3.32} implies that for $i \leq j$
\begin{align*}
    &\nabla_{A_0, \partial_{y_{j+1}}} (\nabla_{A_0, \partial_{y_i}}u'_{j+1}) = (\nabla_{A_0, \partial_{y_i}}u'_{j+1}) \sigma^{(j)}_{j+1} \\
    &\quad \nabla_{A_0, \partial_{y_i}}u'_{j+1}(y_1, ...,y_j, 0, y_{j+2}, ..., y_m) = 0. 
\end{align*}
Uniqueness of solutions to ODE implies that $\nabla_{A_0, \partial_{y_i}}u'_{j+1} = 0$ for $\forall i \leq j$, which further tells us that 
\[
\sigma_i^{(j+1)} = 0 \quad \forall i \leq j+1.
\]
This completes the induction step. 

Due to the uniqueness, the gauge transformation $u$ constructed on $\nu(\gamma) \times X$ is independent of the choices of the identification $\nu(\gamma) \simeq [0, 1] \times I^{m-1}_{\epsilon}$. To see why the restriction $u_y \in \Aut(P_y)$ is independent of the choices of the path $\gamma$, we appeal to the fact that $\Omega$ is simply-connected. Indeed, given any two paths $\gamma_0$ and $\gamma_1$ connecting $y_0$ to $y$, we can find a homotopy of paths $\gamma_t$, $t \in [0,1]$, from $\gamma_0$ to $\gamma_1$ satisfying $\gamma_t(0) = y_0$ and $\gamma_t(1) = y$. Then for any sufficiently close $t_1, t_2 \in [0, 1]$, we can arrange that $\nu(\gamma_{t_1}) \cap \nu(\gamma_{t_2})$ is simply-connected. Thus the uniqueness argument tells us that the gauge transformation $u_t$ constructed on $\nu(\gamma_t) \times X$ has to restrict as the same gauge transformation on $P_y$. 

The gauge transformation constructed above now satisfies $u|_{P_{y_0}} = \id$ and $u\cdot A - A_0 = \tilde{b}_y$. The last condition $d_{H, A_0} \tilde{b}_y = 0$ follows from \eqref{e3.29}. This completes the proof. 
\end{proof}

\autoref{p3.8} enables us to determine when two connections with vanishing $(2,0)$ and $(1,1)$ curvature components are gauge-equivalent. 

\begin{thm}\label{t3.9}
Let $P \to Y \times X$ be a principal $G$-bundle. Suppose $A_0$ and $A$ are two smooth connections on $P$ satisfying
\[
F^{2,0}_{A_0} = F^{2,0}_A = 0 \qquad F^{1,1}_{A_0} = F^{1,1}_A = 0.
\]
Then $u\cdot A = A_0$ for some $u \in \G^f_P$ if and only if 
\begin{equation*}
(1) \; A|_{P_{y_0}} = A_0|_{P_{y_0}} \qquad (2)\; \mathfrak{H}_A(\gamma) = \mathfrak{H}_{A_0}(\gamma) \quad \forall \gamma \in L(Y, y_0). 
\end{equation*}
Moreover, if such a fibre-based gauge transformation $u$ exists, it is unique. 
\end{thm}

\begin{proof}
Recall that a gauge transformation $u \in \G_P$ is called fibre-based if $u|_{P_{y_0}} = \id$. Suppose $u \cdot A = A_0$. Since $u|_{P_{y_0}} = \id$, we know $A|_{P_{y_0}} = A_0|_{P_{y_0}}$. After choosing an identification $P_{y_0} \simeq P_X$, the horizontal holonomy descends to a representation $\pi_1(Y, y_0) \to \G_{P_X}$. Gauge transformations change the representation by the conjugation of their restrictions on $P_{y_0}$. Thus $u|_{P_{y_0}} = \id$ implies that $\mathfrak{H}_A(\gamma) = \mathfrak{H}_{A_0}(\gamma)$. 

Conversely suppose (1) and (2) hold. We construct the gauge transformation $u$ as follows. Given $y_1 \in Y$, we pick an arbitrary smooth path $\alpha$ connecting $y_0$ to $y_1$. With respect to a tubular neighborhood $\nu(\alpha)$ of $\alpha$, \autoref{p3.8} provides us with a unique fibre-based gauge transformation $u^{\alpha}$ defined over $\nu(\alpha) \times X$ so that 
\[
u^{\alpha} \cdot A|_{\nu(\alpha) \times X} - A_0|_{\nu(\alpha) \times X} = b^{\alpha}_y \quad \text{ and } \quad d_{H, A_0} b^{\alpha}_y = 0. 
\]
From the latter relation and the assumption that $b^{\alpha}_{y_0} = 0$, it follows that $b^{\alpha}_y = 0$ for $\forall y \in \nu(\alpha)$.

We claim that the restriction $u^{\alpha}_{y_1}$ of $u^{\alpha}$ at $P_{y_1}$ is independent of the choice of $\alpha$. Suppose we have two paths $\alpha_1$ and $\alpha_2$ connecting $y_0$ to $y_1$. We form a loop $\gamma_o \in L(Y, y_0)$ by concatenating the orientation-reversed path of $\alpha_2$ with $\alpha_1$, and choose a framing of $\gamma_o$ by identifying $\nu(\gamma_o)$ with $[-1, 1] \times I^{m-1}_{\epsilon}$ so that 
\[
\alpha_1(t) = (t-1, 0, ..., 0) \qquad \alpha_2(t) = (1-t, 0, ..., 0) \quad \forall t \in [0, 1]. 
\]
Denote by $u^1$ the gauge transformation constructed on $\nu(\alpha_1) = [-1, 0] \times I^{m-1}_{\epsilon}$ and by $u^2$ the gauge transformation constructed on $\nu(\alpha_2) = [0, 1] \times I^{m-1}_{\epsilon}$. Over $\nu(\gamma_o)$, we write 
\[
A - A_0 = b_y + \sum_{i=1}^m \sigma_idy_i
\]
as before. Note that when we construct $u^i$ over $\nu(\alpha_i)$ in the proof of \autoref{p3.8}, the restriction to $y_1 = (1, 0, ..., 0)$ is given as the initial value $\id$ except at the first step of the induction, where we need to solve the ODE $\nabla_{A_0, \partial_{y_1}} u = u\sigma_1$ with respect to the initial condition $u|_{P_{y_0}} = \id$. By solving this ODE along the concatenated path $\gamma_o$, we conclude that 
\begin{equation}
(u^2_{y_1})^{-1}u^1_{y_1} = \exp\left( \int_{-1}^1 \sigma_1(\gamma_o(t))dt \right) =\mathfrak{H}_{A_0}(\gamma_o)^{-1} \cdot \mathfrak{H}_A(\gamma_o) = \id. 
\end{equation}
This shows that $u^2_{y_1} = u^1_{y_1}$, which establishes the existence of the fibre-based gauge transformation. 

The uniqueness follows directly from \autoref{p3.8} since now we have $u \cdot A = A_0$.
\end{proof}

\subsection{Parallelism}

With the above useful observations in hand, we proceed to prove \autoref{t1.6}. Given a connection $\underline{A}$ on the family bundle $\underline{P}$, we can parallel transport connections on the fibres $P_y$ in the following sense. First, given a path $\alpha$ connecting $y_0$ to $y_1$ in $Y$, the parallel transport map is an bundle isomorphism 
\begin{equation}
\mathfrak{P}_{\underline{A}}(\alpha): P_{y_0} \longrightarrow P_{y_1}. 
\end{equation}
Now, given a connection $B_{y_0}$ on $P_{y_0}$, the horizontal distribution $H_{B_{y_0}} \subseteq TP_{y_0}$ can be parallelly transported to an horizontal distribution in $TP_{y_1}$ by $D\mathfrak{P}_{\underline{A}}(\alpha)(H_{B_{y_0}})$. Suppose further that the holonomy group $\Hol_{\underline{A}} \subseteq \Stab(B_{y_0})$. The parallel transport procedure is then independent of the choice of the paths because
\begin{equation}
D\mathfrak{P}_{\underline{A}}(\alpha_2)^{-1} \circ D\mathfrak{P}_{\underline{A}}(\alpha_1) = \Hol_{\underline{A}}(\alpha_1 * \bar{\alpha}_2) = \id|_{H_{B_{y_0}}}
\end{equation}
for any two paths $\alpha_1$, $\alpha_2$ connecting $y_0$ to $y_1$, where $\alpha_1 * \bar{\alpha}_2$ is the loop obtained by concatenating $\alpha_1$ with the orientation-reversed path of $\alpha_2$. Let us write $\mathfrak{P}^y_{y_0}(H_{B_{y_0}})$ for the horizontal distribution in $P_y$ obtained from parallelly transporting $B_{y_0}$ in the manner described. Then fibrewise summing up the horizontal distributions $H_{\underline{A}}|_{P_y} \oplus \mathfrak{P}^y_{y_0}(H_{B_{y_0}})$ defines an $G$-invariant horizontal distribution in $P$. 

\begin{proof}[Proof of \autoref{t1.6}]
Let $A$ be a smooth $G$-connection on $P$ satisfying 
\[
F^{2, 0}_A = 0 \qquad F^{1, 1}_A = 0.
\]

From $F^{2,0}_A = 0$ it follows by \autoref{l3.1} that the induced connection $\underline{A}$ on the family bundle $\underline{P}$ is flat. As such, its holonomy defines a representation $\rho=\mathfrak{H}_A: \pi_1(Y, y_0) \to \G_{P_X}$. It is then a standard fact (cf. \cite[Chapter 13]{T11}) that such a flat fibre bundle arises from the Borel construction $\tilde{Y} \times_{\rho} P_X$. 

Denote by $B = A|_{P_{y_0}}$ the restricted connection. \autoref{p3.7} tells us that $\im \rho = \Hol_{\underline{A}} \subseteq \Stab(B)$. Then parallelly transporting $B$ with respect to $\underline{A}$ yields a connection $A'$ on $P$. From the construction, we know $\underline{A}' = \underline{A}$, $F^{1,1}_{A'} = 0$, and $A'|_{P_{y_0}} = A|_{P_{y_0}} = B$. Then \autoref{t3.9} tells us there exists a unique fibre-based gauge transformation $u$ on $P$ so that $u \cdot A = A'$. 
\end{proof}

\section{Dimension Reduction Phenomena}\label{s4}

In this section, we first prove vanishing results for the moduli space of $\theta$-instantons over product manifolds, which justify the necessity for imposing the assumptions in \autoref{t1.1}. Next, we apply techniques developed in \autoref{s3} to prove \autoref{t1.1} by relating the $\alpha$-null condition to corresponding curvature conditions. Finally, we prove the compactness result \autoref{t1.4}. 

\subsection{Vanishing of Moduli Spaces}

Let $M = Y \times X$ be the product of two closed Riemannian manifolds with $\dim X \geq 4$. Suppose $\theta \in \Omega^{n-4}(M)$ is a codimension-$4$ parallel form on $M$ given by 
\begin{equation}
\theta = \vol_Y \wedge \beta + \alpha, 
\end{equation}
where $\beta$ is a codimension-$4$ form on $X$, and $\alpha$ does not involve the top form of the $Y$-factor. Since $\theta$ is parallel by assumption, we know that $\beta$ and $\alpha$ are both parallel forms. Let us denote by $\eta_{\beta}$ the smallest eigenvalue of the pointwise defined operator on $X$
\begin{equation}
\begin{split}
\mathcal{Q}_{\beta}: \Omega^2(X) & \longrightarrow \Omega^2(X) \\
\mu & \longmapsto \star_X( \beta \wedge \mu).
\end{split}
\end{equation}
Recall that $\lambda_{\theta}$ is the smallest eigenvalue of $\mathcal{Q}_{\theta}$. We have the following simple observation. 

\begin{lem}\label{l4.1}
Under the above assumptions, we have $\eta_{\beta} \geq \lambda_{\theta}$. 
\end{lem}

\begin{proof}
Let $\mu \in \Omega^2(X)$ be a $2$-form on $X$ that satisfies $\mathcal{Q}_{\beta}(\mu) = \eta_{\beta} \cdot \mu$ and whose pointwise norm equals $1$. Then we know that 
\begin{equation}
\mu \wedge \mu \wedge \beta = \mu \wedge \star_X \mathcal{Q}_{\beta}(\mu) = \eta_{\beta} \cdot \vol_X. 
\end{equation}
Denote by $\omega = \pi_X^*\mu \in \Omega^{0, 2}(M)$ the pull-back form on $M$, which has unit pointwise norm as well. Note that 
\begin{equation}
\omega \wedge \omega \wedge \theta = \omega \wedge \star \mathcal{Q}_{\theta}(\omega) \geq \lambda_{\theta} \cdot \vol_M. 
\end{equation}
On the other hand, $\omega \wedge \omega \wedge \alpha = 0$ (due to degree reasons), which implies that 
\begin{equation}
\omega \wedge \omega \wedge \theta = \pi_X^*\mu \wedge \pi_X^*\mu \wedge \vol_Y \wedge \beta = \eta_{\beta} \cdot \vol_M. 
\end{equation}
We conclude that $\eta_{\beta} \geq \lambda_{\theta}$. 
\end{proof}

Due to \autoref{l4.1}, we may normalize $\theta$, multiplying by a real factor, so that 
\begin{equation}
0 > \eta_{\beta} \geq \lambda_{\theta} = -1. 
\end{equation}
Denote by $\eta_{\beta} = \eta_0 < \eta_1 < ... < \eta_l$ the eigenvalues of $\mathcal{Q}_{\beta}$. Then for $\mu \in \Omega^2(X)$, we have the following inequality
\begin{equation}
\begin{split}
-\eta_{\beta} \cdot |\pi_{\beta}(\mu)|^2 \vol_X & = \sum_{i=1}^l (\eta_i - \eta_{\beta})|\pi_{\eta_i}(\mu)|^2 \cdot \vol_X - \mu \wedge \mu \wedge \beta \\
& \geq - \mu \wedge \mu \wedge \beta. 
\end{split}
\end{equation}
Since $\eta_{\beta} < 0$, we conclude that 
\begin{equation}\label{e4.8}
|\pi_{\beta}(\mu)|^2 \vol_X \geq \frac{\mu \wedge \mu \wedge \beta}{\eta_{\beta}}. 
\end{equation}

Let $\pi: P \to M$ be a principal $SU(r)$-bundle. We write $\pi_X: P_X \to X$ for a reference bundle that is isomorphic to the restriction of $P$ to the fibres $X_y$. Then the $\theta$-instanton charge of $P$ splits into the sum
\begin{equation}\label{e4.9}
\kappa_{\theta}(P) = \Vol(Y) \cdot \kappa_{\beta}(P_X) + \kappa_{\alpha}(P). 
\end{equation}
Appealing to \eqref{e4.8}, we prove the following vanishing result for the moduli space of $\theta$-instantons. 

\begin{prop}\label{p4.2}
With the set-up above, if we further assume that 
\begin{enumerate}
\item either $P$ is $\alpha$-null and $0>\eta_{\beta} > \lambda_{\theta} = -1$;
\item or $P$ is $\alpha$-positive, $0>\eta_{\beta} > \lambda_{\theta} = -1$, and $\kappa_\alpha(P)<\Vol(Y)\cdot \kappa_\beta\cdot (-1-1/\eta_\beta)$;
\item or $P$ is $\alpha$-negative, 
\end{enumerate}
then the moduli space of $\theta$-instantons on $P$ is empty, i.e. $\M_{\theta}(M, P) = \varnothing$. 
\end{prop}

\begin{proof}
(a) Suppose $P$ is $\alpha$-null and $\eta_{\beta} > \lambda_{\theta}$. Let $A$ be a $\theta$-instanton on $P$. Then \autoref{p2.3} and \eqref{e4.9} tell us that 
\begin{equation}
\int_M |F_A|^2 \vol_M = 16r\pi^2\langle c_2(P) \smile [\theta], [M] \rangle  = 16r\pi^2 \Vol(Y) \cdot \kappa_{\beta}(P_X). 
\end{equation}
Let us write $B_y = A|_{P_y}$. Then \eqref{e4.8} tells us that 
\begin{equation}
\begin{split}
\int_M |F_A|^2 \vol_M & \geq \int_M |F_A^{0,2}|^2 \vol_M \\
& = \int_Y \left(\int_{X_y} |\pi_{\beta}(F_{B_y})|^2 \vol_X \right) \vol_Y\\
& \geq -\frac{16r\pi^2}{\eta_{\beta}} \Vol(Y) \cdot \kappa_{\beta}.
\end{split}
\end{equation}
This inequality leads to a contradiction if $0 > \eta_{\beta} > \lambda_{\theta} = -1$.

(b) Suppose that $P$ is $\alpha$-positive and that $0>\eta_{\beta} > \lambda_{\theta} = -1$. By the same discussion as in part (a), we obtain
\begin{equation*}
    \int_M |F_A|^2 \vol_M= 16r\pi^2 (\Vol(Y) \cdot \kappa_{\beta}+\kappa_\alpha)\geq -\frac{16r\pi^2}{\eta_{\beta}} \Vol(Y) \cdot \kappa_{\beta},
\end{equation*}
which gives a contradiction if $\kappa_\alpha(P)<\Vol(Y)\cdot \kappa_\beta\cdot (-1-1/\eta_\beta)$.

(c) Suppose $P$ is $\alpha$-negative. Let $A$ be a $\theta$-instanton on $P$. Then \autoref{p2.3} and \eqref{e4.9} tell us that 
\begin{equation}
\int_M |F_A|^2 \vol_M = 16r\pi^2\langle c_2(P) \smile [\theta], [M] \rangle  < 16r\pi^2 \Vol(Y) \cdot \kappa_{\beta}(P_X). 
\end{equation}
On the other hand, the fact that $0 > \eta_{\beta} \geq \lambda_{\theta} = -1$ implies that 
\begin{equation}
\int_M |F_A|^2 \vol_M \geq -\frac{16r\pi^2}{\eta_{\beta}} \Vol(Y) \cdot \kappa_{\beta} \geq 16r\pi^2 \Vol(Y) \cdot \kappa_{\beta},
\end{equation}
which leads to a contradiction as well. 
\end{proof}

\subsection{Dimension Reduction of Moduli Spaces}

The vanishing result \autoref{p4.2} leads us to consider $\alpha$-null bundles $P$ over products $M = Y \times X$ equipped with the parallel form $\theta = \vol_Y \wedge \beta + \alpha$ satisfying $\eta_{\beta} = \lambda_{\theta} = -1$. That is, $\alpha$-null bundles $P$ over $(M,\theta)$ where $\theta$ is both admissible and compatible with the product structure on $M$. The following observation relates $\theta$-instantons with the curvature conditions in \autoref{t1.6}.

\begin{lem}\label{l4.3}
Let $\pi: P \to M$ be an $\alpha$-null bundle, and suppose that $\theta = \vol_Y \wedge \beta + \alpha$ satisfies $\eta_{\beta} = \lambda_{\theta} = -1$. If $A$ is a $\theta$-instanton on $P$, then
\[
F^{2,0}_A = 0 \quad \text{ and } \quad F^{1,1}_A = 0.
\]
Moreover the restriction of $B_y=A|_{P_y}$ is a $\beta$-instanton for each $y \in Y$. 
\end{lem}

\begin{proof}
We write 
\begin{equation}
\int_M |F_A|^2 \vol_M = \int_M |F_A^{2,0}|^2 \vol_M + \int_M |F_A^{1,1}|^2 \vol_M + \int_M |F^{0,2}_A|^2 \vol_M. 
\end{equation}
The $\alpha$-null assumption implies that 
\begin{equation}
\int_M |F_A|^2 \vol_M = 16r\pi^2 \Vol(Y) \cdot \kappa_{\beta}(P_X).
\end{equation}
The assumption that $\eta_{\beta} = -1$ implies that 
\begin{equation}\label{e4.16}
\int_M |F_A^{0,2}|^2 \vol_M \geq 16r\pi^2 \Vol(Y) \cdot \kappa_{\beta}(P_X). 
\end{equation}
Comparing these three formulae, we conclude that $F^{2,0}_A = F^{1,1}_A = 0$. Moreover, equality holds in \eqref{e4.16}, and so it follows from \autoref{p2.3} that $B_y$ must be a $\beta$-instanton for each $y \in Y$. 
\end{proof}

Appealing to \autoref{t1.6}, we obtain the following criterion that justifies if $\M_{\theta}(M, P)$ is non-empty for $\alpha$-null bundles.

\begin{cor}
Let $\pi: P \to M$ be an $\alpha$-null principal $SU(r)$-bundle, and $\theta = \vol_Y \wedge \beta + \alpha$ satisfies $\eta_{\beta} = \lambda_{\theta} = -1$. Then $\M_{\theta}(M, P) \neq \varnothing$ is and only if the family bundle $\underline{P}$ takes the form $\underline{P} = \tilde{Y} \times_{\rho} P_X$ for some representation $\rho: \pi_1(Y) \to SU(r)$ that stabilizes some $\beta$-instanton on $P_X$.
\end{cor}

It is then natural to ask if the principal bundle $P$ arising this way is topologically isomorphic to the pull-back bundle $\pi_X^*P_X$. From the perspective of characteristic classes, the $\alpha$-null condition merely restricts the second Chern class, and there could be $\alpha$-null bundles whose higher Chern classes decomposes into a mixed type in $H^{s, t}(M; \Z)$. Such $\alpha$-null bundles cannot be a pull-back bundle of the form $\pi_X^*P_X$. However, due to the flatness of $\underline{P}$, we can prove that $P$ has to be the pull-back bundle if $\dim Y \leq 4$, as stated in \autoref{p1.5}. 

\begin{proof}[Proof of \autoref{p1.5}]
Given a representation $\rho: \pi_1(Y) \to SU(r)$, we write $Q_{\rho} := \tilde{Y} \times_{\rho} SU(r)$ for the associated flat $SU(r)$-bundle. Then $\underline{P} = Q_{\rho} \times_{\iota} \G_{P_X}$, where $\iota: SU(r) \to \G_{P_X}$ is the inclusion of constant gauge transformations. When $\dim Y \leq 3$, any $SU(r)$-bundle on $Y$ is topologically trivial. Thus the family bundle $\underline{P}$ is isomorphic to the product $Y \times P_X$, which means $P$ is the pull-back bundle $\pi_X^*P_X$. 

When $\dim Y = 4$, the flatness of $Q_{\rho}$ implies that $c_2(Q_{\rho})  = 0 \in H^4(Y; \Z)$, noting that $H^4(Y; \Z)$ is torsion-free. Since $SU(r)$-bundles on a closed $4$-manifold are classified by the second Chern class, it follows that $Q_{\rho}$ is again topologically trivial. This finishes the proof. 
\end{proof}

Things get slightly more complicated when $\dim Y > 4$. Non-pull-back bundles could arise from the following three possibilities, which one cannot exclude in general. The first is that $H^4(Y; \Z)$ is not torsion-free, in which case the flat bundle $Q_{\rho}$ could have torsion $c_2$. The second comes from non-trivial higher Chern classes $c_l$, $l \geq 3$, of the bundle $Q_{\rho}$. As for the third possiblity, we note that the differential form $\alpha$ might be a sum of mixed types. Thus for an $\alpha$-bundle $P$, the second Chern class $c_2(P)$ might have non-trivial mixed components in $H^{s, t}(M; \Z)$ which cancels out after pairing with $\alpha$. 

Now we proceed to construct all $\theta$-instantons over an $\alpha$-null bundle $P$ whose associated family bundle is given by $\underline{P} = \tilde{Y} \times_{\rho} P_X$. The major difficulty in the construction arises from analyzing the possibilities of the stabilizers when the rank $r$ is large. As mentioned in the introduction, this can be avoided by considering moduli spaces of irreducible connections. 

We know from \cite[Lemma 4.2.8]{DK90} that the stabilizer of a connection $A$ on $P$ is given by the centralizer of its holonomy group $\Hol_A$. In particular, the smallest possible stabilizer of connections on an $SU(r)$-bundle is the center $\Z/r$, a copy of the cyclic group of order $r$. Recall that we call connections with stabilizer $\Z/r$ to be irreducible. 

Given an element $\zeta \in H^1(Y; \Z/r)$, we get a flat line bundle $L_{\zeta} \to Y$ equipped with a flat connection $\Gamma_{\zeta}$ whose holonomy is precisely $\zeta$. Since $c_1(L_{\zeta})$ is $r$-torsion, the twisted bundle $\pi_1^*L_{\zeta} \otimes P$ is isomorphic to $P$, where $\pi_1: Y \times X \to Y$ is the projection to the first factor. Given an $SU(r)$-connection $A$ on $P$, we get a twisted connection $\pi_1^*\Gamma_{\zeta} \otimes A$ on $\pi_1^*L_{\zeta} \otimes P$, which corresponds to an $SU(r)$-connection $A_{\zeta}$ on $P$ once an identification $\pi_1^*L_{\zeta} \otimes P \simeq P$ is chosen. Nevertheless, the gauge equivalence class $[A_{\zeta}]$ is well-defined.

Given a representation $\rho: \pi_1(Y, y_0) \to \Z/r$, we denote by $\underline{P}_{\rho} = \tilde{Y} \times_{\rho} P_X$ the family bundle induced from $\rho$, and $P_{\rho} \to Y \times X$ the corresponding principal $SU(r)$-bundle. Given a connection $B$ on $P_X$, we denote by $\pi^*_{\rho}B$ the connection on $P_{\rho}$ obtained by parallel transporting $B$ using the representation $\rho$. 

With such a set-up, we are led to consider the following map: 
\begin{equation}\label{e4.17}
\begin{split}
\Theta^*_{\rho}: H^1(Y; \Z/r) \times \mathscr{B}^*_{P_X} & \longrightarrow \mathscr{B}^*_{P_{\rho}} \\
(\zeta, [B]) & \longmapsto [(\pi_{\rho}^*B)_{\zeta}],
\end{split}
\end{equation}
where the superscript $^*$ means to take the part of the configuration space $\mathscr{B} = \A/\G$ that consists of irreducible connections. Once an identification $\pi_1^*L_{\zeta} \otimes P \simeq P$ is fixed, the $(1, 0)$-component of $(\pi^*_{\rho}B)_{\zeta}$ is independent of the connection $B$, since its restriction on $Y \times \{x\}$ is determined by $\rho$ and $\zeta$. 

Now we prove the main result \autoref{t1.1}. Let us write $\M^*_{\theta}(M, P_{\rho})$ for the moduli space of irreducible $\theta$-instantons on $M$, and $\mathcal{N}^*_{\beta}(X, P_X)$ for the moduli space of irreducible $\beta$-instantons on $X$. \autoref{t1.1} states that the map 
\begin{equation}
\Theta^*_{\rho}: H^1(Y; \Z/r) \times \mathcal{N}^*(X, P_X) \longrightarrow \M^*_{\varphi}(M, P_{\rho})
\end{equation}
induces an isomorphism between moduli spaces. The proof consists of two parts. In the first part, we prove that $\Theta^*_{\rho}$ is a homeomorphism with respect to $C^{\infty}$-topology. In the second part, we prove that its differential identifies the Zariski tangent spaces. 

\begin{proof}[Proof of \autoref{t1.1} (Part I)]
Given $\zeta \in H^1(Y; \Z/r)$ and $B$ an irreducible $\beta$-instanton on $P_X$, we write $A_{\zeta} = (\pi^*_{\rho}B)_{\zeta}$ for the induced connection on $P_{\rho}$ after choosing an identification $\pi_1^*L_{\zeta} \otimes P \simeq P$. Then the construction of $A_{\zeta}$ tells us that 
\[
F^{2,0}_{A_{\zeta}} = 0 \quad \text{ and } \quad F^{1,1}_{A_{\zeta}} = 0.
\]
Since $A_{\zeta}|_{P_y}$ is identified with $B$, we see that $A_{\zeta}$ minimizes the Yang--Mills functional on $P_{\rho}$. Thus $A_{\zeta}$ is an irreducible $\theta$-instanton on $P_{\rho}$. 

Now we fix a basepoint $y_0 \in Y$, and identifies $P_{y_0}$ with $P_X$. 

\vspace{2mm}

{\noindent \bf \em Step 1.(Injectivity)} \hspace{1mm} Given $\zeta \in H^1(Y; \Z/r)$, we write $\rho_{\zeta}: \pi_1(Y, y_0) \to \Z/r$ for the representation obtained by twisting $\rho$ with $\zeta$. Then for any irreducible connection $B$ on $P_X$, the induced connection $A_{\zeta}$ on $P_{\rho}$ satisfies $\Hol^h_{A_{\zeta}} = \rho_{\zeta}$. Thus if $\Theta_{\rho}(\zeta, [B]) = \Theta_{\rho}(\zeta', [B'])$, we know that $\zeta = \zeta'$. Moreover, $B = A_{\zeta}|_{P_{y_0}}$ and $B' = A'_{\zeta'}|_{P_{y_0}}$. So $[A_{\zeta}] = [A'_{\zeta'}]$ implies that $B$ and $B'$ are gauge-equivalent. This proves injectivity.

\vspace{2mm}
{\noindent \bf \em Step 2.(Surjectivity)} \hspace{1mm} Given $[A] \in \M^*_{\theta}(M, P_{\rho})$, we write $B = A|_{P_{y_0}}$. \autoref{t1.6} implies that $\Stab(B) = \Stab(A)$. Thus $B$ is an irreducible $\beta$-instanton. It follows from \autoref{p3.7} that the horizontal holonomy $\Hol^h_A = \Z/r$. Then the difference between the horizontal holonomy $\Hol^h_{A}$ and the representation $\rho$ defines a cohomology class $\zeta \in H^1(Y; \Z/r)$. From the construction, we know that $\Hol^h_{A_{\zeta}} = \Hol^h_A$, and $A|_{P_{y_0}} = A_{\zeta}|_{P_{y_0}} = B$. Then \autoref{t3.9} implies that $[A]  = [A_{\zeta}]$, which proves surjectivity. 

\vspace{2mm}
{\noindent \bf \em Step 3.(Continuity of $\Theta^*_{\rho}$)} \hspace{1mm} Let $\{B_n\}_{n \geq 1}$ be a sequence of irreducible $\beta$-instantons on $P_X$ that converges to $\beta$-instanton $B_o$ in $C^{\infty}$-topology. Since $H^1(Y; \Z/r)$ is finite, we may fix $\zeta \in H^1(Y; \Z/r)$ and write $A_n = (\pi^*_{\rho}B_n)_{\zeta}$. It suffices to prove that $A_n \to A_o$ in $C^{\infty}$-topology as $n \to \infty$. 

Let us write $b_n = B_n - B_o$ and $a_n = A_n - A_o$. Since the induced connections $\underline{A}_n = \underline{A}_o$ agree on the family bundle $\underline{P}_{\rho}$, we see that $a_n \in \Omega^{0, 1}$. Thus we can regard $a_n$ as a section taking value in $\Omega^1(X_y, \g_{P_y})$ at $y$. Then $a_n$ is characterized by
\begin{equation}
\nabla_{\underline{A}_o} a_n = 0 \qquad a_n|_{P_{y_0}} = b_n. 
\end{equation}
Since $Y$ is compact, we conclude that for all $k \in \N$ one can find a constant $c_k > 0$ such that 
\begin{equation}\label{e4.20}
\left\| a_n \right\|_{L^2_k(Y \times X)} \leq c_k \left\|b_n \right\|_{L^2_k(X)},
\end{equation}
where the Sobolev norms on sections of $\g_{\underline{P}_{\rho}}$ and $\g_{P_X}$ are defined using $A_o$ and $B_o$ respectively. Then it follows from the Sobolev embedding theorem that $a_n \to 0$ in $C^{\infty}$-topology as $n \to \infty$. 

\vspace{2mm}
{\noindent \bf \em Step 4.(Continuity of $(\Theta^*_{\rho})^{-1}$)} \hspace{1mm} Let $\{A_n\}_{n\geq 1}$ be a sequence of irreducible $\theta$-instantons on $P_{\rho}$ that converges to some irreducible $\theta$-instanton $A_o$ in $C^{\infty}$-topology. We write $B_n = A_n|_{P_{y_0}}$. Then it clear that $B_n \to B_o$ in $C^{\infty}$-topology. The element $\zeta_n \in H^1(Y; \Z/r)$ is the difference between $\Hol^h_{A_n}$ and $\rho$. Since $H^1(Y; \Z/r)$ is finite, convergence of $A_n$ implies that $\zeta_n = \zeta_o$ when $n$ is sufficiently large. 
\end{proof}

To establish the second part that $\Theta^*_{\rho}$ identifies the Zariski tangent spaces, we start with two auxiliary lemmas.

\begin{lem}\label{l4.4}
Let $A$ be a connection on $P \to Y \times X$, and $B_y = A|_{P_y}$. Given $a \in \Omega^1(M, \g_P)$, we have the implication
\[
d^*_A a = 0 \;\;\Longrightarrow\;\; d^*_{B_y} b_y = 0, \; \forall y \in Y,
\] 
where $b_y = a|_{X_y}$.
\end{lem}

\begin{proof}
Let us write $a = a^{1, 0} + a^{0, 1}$. Then $\star_{\varphi} a = \star_Y a^{1,0} \wedge \vol_X - \vol_Y \wedge \star_X a^{0, 1}$. From the decomposition $d_A = d_{H, A} + d_{f, A}$, it follows that 
\[
d_{H, A} \star_Y a^{1,0} = 0 \text{ and } d_{f, A} \star_X a^{0, 1} = 0
\]
whenever $d \star_{\varphi} a = 0$. The result now follows from $d_{B_y} \star_X b_y = (d_{f, A} \star_X a^{0, 1})|_{X_y} = 0$.
\end{proof}

Recall that the deformation complex of $\M_{\theta}(M, P)$ at a $\theta$-instanton is given by 
\begin{equation}\tag{$E_{\theta, A}(M)$}
0 \to \Omega^0(M, \g_P) \xrightarrow{-d_A} \Omega^1(M, \g_P) \xrightarrow{\pi_{\theta^{\perp}} \circ d_A} \Omega^2(M, \g_P) \to 0. 
\end{equation}
Replacing $A$ by a $\beta$-instanton $B$, and $M$ by the fibre $X$, we get the deformation complex $E_{\beta, B}(X)$ of $\mathcal{N}_{\beta}(X, P_X)$ at $[B]$. The next observation is more crucial and makes use of the $\alpha$-null condition on $P$ in an essential way. 

\begin{lem}\label{l4.5}
Let $(M, \theta)$ be the pair as above, and $A$ be a $\theta$-instanton on an $\alpha$-null $SU(r)$-bundle $P$. Suppose $a \in \Omega^1(M, \g_P)$ satisfies $\pi_{\theta^{\perp}}(d_A a) = 0$. Then 
\[
(d_A a)^{2,0} = 0 \qquad (d_A a)^{1,1} = 0 \qquad \pi_{\beta^{\perp}}(d_{B_y} b_y) = 0, \; \forall y \in Y,
\]
where $B_y = A|_{P_y}$, and $b_y = a|_{X_y}$. 
\end{lem}

\begin{proof}
Considering the curvature $F_{A + a}$, \eqref{e2.5} tells us that 
\begin{equation}
-2r\tr(F_{A+a} \wedge F_{A+a}) \wedge \theta = \left(\sum_{i=1}^k \lambda_i \cdot |\pi_{\lambda_i}(F_{A+a})|^2 - |\pi_{\theta}(F_{A+a})|^2 \right) \vol_M.
\end{equation}
Thus the Yang--Mills functional evaluated at the connection $A+a$ is 
\begin{equation}\label{e4.22}
\|F_{A+a}\|^2_{L^2(M)} = \sum_{i=1}^k (\lambda_i + 1)\|\pi_{\lambda_i}(F_{A+a})\|^2_{L^2(M)} + 16r\pi^2 \kappa_{\theta}(P).
\end{equation}
Since $A$ is a $\theta$-instanton, we know that $\pi_{\theta^{\perp}}(F_A) = 0$. Thus we obtain 
\begin{equation}\label{e4.23}
\pi_{\theta^{\perp}}(F_{A+a}) = \pi_{\theta^{\perp}}(F_A + d_A a + a\wedge a)  = \pi_{\theta^{\perp}}(a \wedge a).
\end{equation}
Appealing to \eqref{e2.5} again, we have 
\begin{equation}\label{e4.24}
 0 = - 2r\tr\left((a \wedge a) \wedge (a \wedge a) \right) = \left(\sum_{i=1}^k \lambda_i \cdot |\pi_{\lambda_i}(a \wedge a)|^2 - |\pi_{\theta}(a \wedge a)|^2 \right) \vol_M.
\end{equation}
It follows from \eqref{e4.22}---\eqref{e4.24} that 
\begin{equation}\label{e4.25}
\|F_{A+a}\|^2_{L^2(M)} = \| a \wedge a\|^2_{L^2(M)} + 16r\pi^2 \kappa_{\theta}(P). 
\end{equation}
Repeat the derivation of \eqref{e4.22} on $B_y + b_y$, we see that 
\begin{equation}\label{e4.26}
\|F_{B_y+b_y}\|^2_{L^2(X_y)} = \sum_{j=1}^l(\eta_j + 1) \|\pi_{\eta_j}(F_{B_y+b_y})\|^2_{L^2(X_y)} + 16r\pi^2\kappa_{\beta}(P_X). 
\end{equation}
The $\alpha$-null assumption implies that $\kappa_{\theta}(P) = \Vol(Y) \cdot \kappa_{\beta}(P_X)$. Thus it follows from \eqref{e4.25} and \eqref{e4.26} that 
\begin{equation}\label{e4.27}
\begin{split}
\| a \wedge a\|^2_{L^2(M)} & = \|F^{2,0}_{A+a}\|^2_{L^2(M)} + \|F^{1, 1}_{A+a}\|^2_{L^2(M)}  \\
& + \sum_{j=1}^l \int_Y (\eta_j + 1)\|\pi_{\eta_j}(F_{B_y+b_y})\|^2_{L^2(X_y)} \vol_Y. 
\end{split}
\end{equation}
Since $F^{2, 0}_A = F^{1, 1}_A = 0$, and $\pi_{\beta^{\perp}}(F_{B_y}) = 0$, the right-hand side of \eqref{e4.27} can be further reduced to 
\begin{equation}\label{e4.28}
\begin{split}
\| a \wedge a\|^2_{L^2(M)} & = \| (d_Aa + a\wedge a)^{2, 0}\|^2_{L^2}(M) + \| (d_A a + a\wedge a)^{1,1}\|^2_{L^2(M)} \\
& + \sum_{j=1}^l \int_Y (\eta_j + 1)\|\pi_{\eta_j}(d_{B_y}b_y + b_y \wedge b_y)\|^2_{L^2(X_y)} \vol_Y. 
\end{split}
\end{equation}
Replacing $a$ in \eqref{e4.28} by $\epsilon a$ for small positive constant $\epsilon$, we get 
\begin{equation}
\begin{split}
& \epsilon^2 \left( \|(d_A a)^{2,0}\|^2_{L^2(M)} + \|(d_A a)^{1, 1}\|^2_{L^2(M)}\right) \\
& + \epsilon^2 \sum_{j=1}^l \int_Y (\eta_j + 1)\|\pi_{\eta_j}(d_{B_y}b_y)\|^2_{L^2(X_y)} \vol_Y = \mathcal{O}(\epsilon^3).
\end{split}
\end{equation}
Dividing through by $\epsilon^2$ and taking the limit as $\epsilon\rightarrow 0$, it follows that 
\begin{equation}
(d_A a)^{2, 0} = 0 \qquad (d_A a)^{1, 1} = 0 \qquad \pi_{\beta^{\perp}}(d_{B_y} b_y) = 0
\end{equation}
as desired. 
\end{proof}

Now we are ready to prove the second part of \autoref{t1.1}, which amounts to say that the differential of $D\Theta^*_{\rho}$ defines an isomorphism 
\begin{equation}
D\Theta^*_{\rho}|_{(\zeta, [B])}: H^1(E_{\beta, B}(X)) \longrightarrow H^1(E_{\theta, A}(M)),
\end{equation}
where $[A] = \Theta^*_{\rho}(\zeta, [B])$. 

\begin{proof}[Proof of \autoref{t1.1} (Part II)]
The first cohomology $H^1(E_{\theta, A}(M))$ is identified as $\ker (\pi_{\theta^{\perp}} \circ d_A) \cap \ker (-d^*_A)$, so is $H^1(E_{\beta, B}(X))$. Denote by 
\[
\Pi_2: H^1(Y; \Z/r) \times \mathcal{N}^*_{\beta}(X, P_X) \to \mathcal{N}^*(X, P_X)
\]
the projection map to the second factor. The composition $\Pi_2 \circ (\Theta^*_{\rho})^{-1}$ is the restriction map $\mathfrak{R}: [A] \mapsto [A|_{P_{y_0}}]$. Then the differential $D\mathfrak{R}|_{[A]}$ is again the restriction map $a \mapsto a|_{P_{y_0}}$ under the above identifications. 

\autoref{l4.4} and \autoref{l4.5} tell us that $D\mathfrak{R}|_{[A]}$ restricts to a linear map from $H^1(E_{\theta, A}(M))$ to $H^1(E_{\beta, B}(X))$. Since $H^1(Y; \Z/r)$ is discrete, the differential of $\Theta^*_{\rho}$ is the inverse of $D\mathfrak{R}$, which then must be an isomorphism between the cohomology space. 
\end{proof}

Although it is quite complicated to exhaust all possibilities of stabilizers, we can still determine reducible $\theta$-instantons when specializing to $SU(2)$-bundles. 

Let $P \to M$ be a principal $SU(2)$-bundle and $E = P \times_{SU(2)} \C^2$ be the associated vector bundle. Recall a connection on $P$ is called abelian if its stabilizer is $U(1)$. An abelian connection $A$ on $P$ decomposes the vector bundle into the sum of two line bundles
\[
E = L \oplus \bar{L}, \text{ with } c_2(E) = -c_1^2(L).
\]
The connection decomposes accordingly as $A = A^{\dagger} \oplus \bar{A}^{\dagger}$. We write 
\begin{equation}
\Pic(Y) = H^1(Y; \R) / H^1(Y; \Z)
\end{equation}
for the Picard torus of $Y$, which parametrizes flat $U(1)$-connections on the trivial line bundle over $Y$. Given an element $\xi \in \Pic(Y)$, we choose a flat connection $\Gamma_{\xi}$ representing $\xi$ and write $A^{\dagger}_{\xi} = \pi_1^*\Gamma_{\xi} \otimes A^{\dagger}$ for the twisted connection on $L \to M$. We write $A_{\xi} = A^{\dagger}_{\xi} \oplus \bar{A}^{\dagger}_{\xi}$ for the induced $SU(2)$-connection on $P$. Then the gauge equivalence class $[A_{\xi}]$ is independent of the choices of the flat connection $\Gamma_{\xi}$ representing $\xi$. 

To define the twisting map $\Theta$ in the abelian case, we form an $SU(2)$-bundle $P_{\rho} \to M$ with respect to a representation $\rho: \pi_1(Y, y_0) \to U(1) \subseteq SU(2)$ and a principal $SU(2)$-bundle $P_X \to X$. The parallel transport of a reducible connection $B$ on $P_X$ results in a reducible connection on $P_{\rho}$, which we denote by $\pi^*_{\rho}B$. Then we define 
\begin{equation}
\begin{split}
\Theta^a_{\rho}: \Pic(Y) \times \mathscr{B}^a_{P_X} & \longrightarrow \mathscr{B}^a_{P_{\rho}} \\
(\xi, [B]) & \longmapsto [(\pi^*_{\rho}B)_{\xi}],
\end{split}
\end{equation}
where the superscript $^a$ means restriction to the part of the configuration space that consists of abelian connections. Replacing the finiteness of $H^1(Y; \Z/r)$ with the compactness of $\Pic(Y)$, the same argument in the first part proof of \autoref{t1.1} proves the corresponding statement in \autoref{t1.2} for the abelian locus that the map
\[
\Theta^a_{\rho}: \Pic(Y) \times \mathcal{N}^a_{\beta}(X, P_X) \longrightarrow \M^a_{\varphi}(M, P_{\rho})
\]
is a homeomorphism.

To obtain the second part of \autoref{t1.2}, we analyze the effect of the differential of $\Theta^a_{\rho}$ on Zariski tangent spaces of the abelian strata. An abelian $\theta$-instanton $A$ on $P$ splits the adjoint bundle as $\g_P = \underline{\R} \oplus L^{\otimes 2}$. Its restriction $B$ on $P_X$ also splits the adjoint bundle $\g_{P_X} = \underline{\R} \oplus L^{\otimes 2}_X$. We write the induced splittings on the deformation complexes as 
\[
E_{\beta,B}(X) = E_{\beta, B}^{\R}(X) \oplus E_{\beta, B}^{L}(X) \qquad E_{\theta, A}(M) = E_{\theta, A}^{\R}(M) \oplus E_{\theta, A}^L(M). 
\]
The second part proof of \autoref{t1.1} tells us that the differential of $\Theta^a$ at $(\xi, [B])$ identifies $H^1(E^L_{\beta, B}(X))$ with $H^1(E^L_{\theta, A}(X))$. On the other hand, we have the following identifications of $H^1(E_{\beta, B}^{\R}(X))$ and $H^1(E^{\R}_{\theta, A}(M))$.

\begin{lem}\label{l4.6}
With the set-up above, we have
\[
H^1(E_{\beta, B}^{\R}(X)) = H^1(X; \R) \qquad H^1(E_{\theta, A}^{\R}(M)) = H^1(M; \R).
\]
\end{lem}

\begin{proof}
We note that $E^{\R}_{\theta, A}(M)$ is 
\[
0 \to \Omega^0(M, \R) \xrightarrow{-d} \Omega^1(M, \R) \xrightarrow{\pi_{\theta^{\perp}} \circ d} \Omega^2_{\theta^{\perp}}(M, \R) \to 0. 
\]
When $\pi_{\theta^{\perp}}(da) = 0$, \eqref{e2.5} tells us that 
\[
0 = \int_M da \wedge da \wedge \theta = - \int_M |\pi_{\theta}(da)|^2 \vol_M,
\]
from which we see that $\pi_{\theta}(da) = 0$. Since $\pi_{\theta^{\perp}}(da) = 0$ implies that $da=0$, we conclude that  $H^1(E_{\theta, A}^{\R}(M)) = H^1(M; \R)$. Replacing $(M, \theta)$ by $(X, \beta)$, we get the other identification as well. 
\end{proof}

We note that $T_{\xi} \Pic(Y) = H^1(Y; \R)$ for $\forall \xi \in \Pic(Y)$. From the construction of $\Theta^a$, it is straightforward to see that the effect of the differential $D\Theta^a$ on $T_{\xi} \Pic(Y)$ is given by the pull-back map $z \mapsto \pi_1^*z$. Similarly the restriction of $D\Theta^a$ on $H^1(E^{\R}_{\beta, B}(X))$ is also the pull-back map $\beta \mapsto \pi_2^* \beta$. Thus the differential $D\Theta^a|_{(\xi, [B])}$ restricts to the isomorphism
\begin{equation}
\begin{split}
H^1(Y; \R) \oplus H^1(X; \R) & \longrightarrow H^1(M; \R) \\
([z], [\beta]) & \longmapsto [\pi^*_1z + \pi^*_2 \beta]. 
\end{split}
\end{equation}
We summarize the discussion above in the following proposition, which completes the proof of \autoref{t1.2}. 

\begin{prop}\label{p4.7}
With the set-up above, the differential of $\Theta^a_{\rho}$ defines an isomorphism 
\[
D\Theta^a_{\rho}|_{(\xi, [B])}: H^1(Y; \R) \oplus H^1(E_{\beta, B}(X)) \longrightarrow H^1(E_{\theta, A}(M)),
\]
where $[A] = \Theta^a_{\rho}(\xi, [B])$.  
\end{prop}

For $SU(2)$-connections, the only possible stabilizers are the center $\Z/2$, the maximal torus $U(1)$, and the whole group $SU(2)$. We have discussed the first two cases. The last case is straightforward. 

\begin{dfn}
An $SU(2)$-connection $A$ is called central if $\Stab(A) = SU(2)$. 
\end{dfn}

Suppose $A$ is a central $\theta$-instanton on a principal $SU(2)$-bundle $P \to M$. Then $A$ can be used to trivialize the corresponding vector bundle, which implies that $c_2(P) = 0$. It follows from \autoref{p2.3} that $A$ is a flat connection. Thus central $\theta$-instantons are flat connections on the trivial $SU(2)$-bundle with $\Hol(A) \subseteq \Z/2$. Let us write 
\begin{equation}
\chi^c(Y):= \Hom\left(\pi_1(Y, y_0), \Z/2\right)
\end{equation}
for the character variety of $Y$, which parametrizes gauge equivalence classes of flat $SU(2)$-connections on $Y$ with holonomy contained in the center $\Z/2$. Then for the trivial bundle $P \to M$, we get the identification between the moduli spaces of central instantons immediately from the fact that $\pi_1(M) \simeq \pi_1(Y) \times \pi_1(X)$. Thus, we have an isomorphism
\begin{equation}
\Theta^c: \chi^c(Y) \times \mathcal{N}^c_{\beta}(P_X) \longrightarrow \M^c_{\theta}(M, P),
\end{equation}
where the superscript $^c$ means to take the part of the configuration space consisting of central connections. 

In principle, we could run through the discussion above to define the map $\Theta^{\Gamma}_{\rho}$, which identifies the moduli space of instantons with stabilizer $\Gamma$ given any fixed subgroup $\Gamma \subseteq SU(r)$. All one needs to do is to analyze how the $\Gamma$-stabilized connections decompose the corresponding vector bundle, as this provides the information of the flat connections on $Y$ that can be used to twist the instantons on $P_X$. 

\subsection{Compactification of Moduli Spaces}

Now let $X$ be a closed $4$-manifold in the Riemannian product $M = Y \times X$. To ensure that $\eta_{\beta} = -1$, we require that $\beta \equiv 1$ is the constant $1$ function. Then we know $\theta = \vol_Y + \alpha$. Suppose now that $P \to M$ is an $\alpha$-null principal $SU(r)$-bundle with $\M_{\theta}(M, P) \neq \varnothing$. \autoref{t1.6} tells us that $P = P_{\rho}$ for some representation $\rho: \pi_1(Y, y_0) \to SU(r)$. 

Let us write $\langle c_2(P_X), [X] \rangle = \kappa \in \N$. For each integer $0 \leq s \leq \kappa$, we write $P_X(s) \to X$ for the bundle characterized by $\langle c_2(P_X(s)), [X] \rangle  = \kappa - s$, and $P_{\rho}(s)$ the principal $SU(r)$-bundle corresponding to the family $\underline{P}_{\rho}(s) = \R^3 \times_{\rho} P_X(s)$. We note that $P_{\rho}(s)$ is still an $\alpha$-null bundle whose second Chern class is given by 
\begin{equation}
c_2(P_{\rho}) - c_2(P_{\rho}(s)) = s\PD([Y]) \in H^2(Y \times X; \Z). 
\end{equation}
As introduced in \autoref{s1}, we have the moduli space of ideal $\theta$-instantons quipped with the weak topology given by \autoref{d1.3}:
\begin{equation}
\overline{\M}_{\theta}(M, P_{\rho})= \M_{\theta}(X, P_{\rho})\cup\bigcup_{s=1}^{\kappa} \M_{\theta}(M, P_{\rho}(s)) \times \Sym^s(X). 
\end{equation}
by adjoining ideal $\theta$-instantons to the original moduli space. We now prove \autoref{t1.4} that this moduli space is compact. 

\begin{proof}[Proof of \autoref{t1.4}]
Let $\{[A_n]\}_{n \geq 1}$ be a sequence of $\theta$-instantons in $\M_{\theta}(M, P)$. \autoref{t1.6} provides us with a sequence of representations $\rho_n: \pi_1(Y, y_0) \to SU(r)$, and ASD connections $B_n = A_n|_{P_{y_0}}$ such that $A_n$ arises as the parallel transport of $B_n$ using the flat connection corresponding to $\rho_n$. Since the representation variety 
\[
\left\{\rho: \pi_1(Y, y_0) \to SU(r) \right\} / \text{ conj. }
\]
is compact, we may assume that $\rho_n \to \rho_o$ as $n \to \infty$ after passing to a subsequence and taking conjugations. On the other hand, Uhlenbeck's compactness implies that one can find $\pmb{x} \in \Sym^s(X)$ and $[B_o] \in \mathcal{N}(P_X(s))$ so that $[B_n]$ converges to the ideal instanton $([B_o], \pmb{x})$ after further passing to a subsequence. Since the stabilization equation $d_{B_n} u_n = 0$ is a closed condition, we know that the gauge transformation given by the horizontal holonomy of $\rho_o$ stabilizes $B_o$. We denote by $A_o$ the connection obtained by parallel transporting $B_o$ via $\rho_o$, which is a $\theta$-instanton on $P(s)$ by construction. 

Let us write $B_{n, y} = A_n|_{P_y}$. Note that $P$ is $\alpha$-null. Thus $F_{A_n} = F^{2,0}_{A_n}$. Appealing to \autoref{t1.6}, we get $|F_{B_{n, y}}|  = |F_{B_n}|$. Given $f \in C^{\infty}(M)$, we write $f_y = f|_{X_y}$. Then 
\begin{equation}
\begin{split}
\int_M f \cdot |F_{A_n}|^2 \vol_{\varphi} & = \int_Y \left(\int_{X_y} f_y \cdot |F_{B_{n, y}}|^2 \vol_{X_y} \right)\vol_{Y} \\
& = \int_Y \left(\int_{X_y} f_y \cdot |F_{B_n}|^2 \vol_{X_y} \right)\vol_{Y} \\
& \to \int_Y \left(\int_{X_y} f_y \cdot |F_{B_o}|^2 \vol_{X_y} \right) \vol_Y + \sum_{x \in \pmb{x}} \int_Y f \vol_Y \\
& = \int_M f \cdot |F_{A_o}|^2 \vol_{\varphi} + \sum_{x \in \pmb{x}} \int_Y f \vol_Y. 
\end{split}
\end{equation}
On the other hand, given any compact subset $K \subseteq M \backslash Y \times \pmb{x}$, we can choose a compact subset $K_X \subseteq X \backslash \pmb{x}$ so that $K \subseteq Y \times K_X$. After applying gauge transformations, we know that $B_n$ converges to $B_o$ in $C^{\infty}(K_X)$. It follows from the estimate \eqref{e4.20} that $A_n$ converges to $A_o$ in $C^{\infty}(Y \times K_X)$ up to gauge transformations. Thus, we have proven that the closure of $M_{\theta}(X,P)$ lies inside $\overline{M}_{\theta}(X,P)$ is compact. However, by the gluing results of Taubes \cite{T82}, any ideal instanton can occur as a limit. This completes the proof.
\end{proof}

Note that $\theta = \vol_Y + \alpha$. The restriction of $\theta$ to the slice $Y \times \{x\}$ is the volume form on $Y$, which is equivalent to say that $Y \times \{x\}$ is calibrated by $\theta$. Tian \cite{T00} proved that the bubbling set of ideal $\theta$-instantons must be calibrated by $\theta$. \autoref{t1.6} tells us that when $P$ is $\alpha$-null, and $\dim X = 4$, these are the only calibrated submanifolds that can form the bubbling set of ideal instantons. Depending on the context, there might be lots of other calibrated submanifolds, which we shall see in the following section of examples. 

\section{Hyperk\"ahler Fibrations}\label{s5}

The motivating examples of our study on $\theta$-instantons come from hyperk\"ahler fibrations in dimension $7$ and $8$. If we require the total space mainfold is equipped with a $G_2$- or $\Spin(7)$-structure so that the fibres are coassociative or Cayley, then we shall prove in this section that the manifold must be a finite quotient of a Riemannian product manifold with certain contraint on the factors. Thus, essentially we only need to deal with instantons over product manifolds for such cases. 

\subsection{Coassociative Fibrations}

In this subsection, we restrict ourselves to the special kind of $G_2$-manifolds that carry fibration structures with coassociative fibres. One may consult \cite{B10, D17} for more details on the general theory of coassociative fibrations and \cite{KL21} for concrete examples. 

We shall write $(M, \varphi)$ for a complete torsion-free $G_2$-manifold (not necessarily compact), and $\psi = \star_{\varphi} \varphi$ the companion $4$-form. 

\begin{dfn}
A $4$-dimensional submanifold $X \subseteq M$ is called coassociative if $\varphi|_X = 0$. 
\end{dfn}

It is a standard fact proved in \cite{HL82} that the condition $\varphi|_X = 0$ is equivalent to the condition that $\psi$ restricts as a volume form on $X$, i.e. $\psi$ is a calibration on $X$. The deformation theory of coassociative submanifolds is well studied in \cite{M98}. Given a coassociative submanifold $X \subseteq M$, the volume form $\psi|_X$ provides us with an orientation. The normal bundle of $X$ is identified with the bundle of self-dual $2$-forms via 
\[
\begin{split}
NX & \longrightarrow \Lambda^+T^*X \\
v & \longmapsto \iota_v \varphi|_X,
\end{split}
\]
where we have identified the normal bundle with the orthogonal complement $(TX)^{\perp_{\varphi}} \subseteq TM|_X$ using the metric $g_{\varphi}$. The main result of McLean \cite{M98} states that the local coassociative deformation moduli near $X$ can be identified with the space $\mathcal{H}^+(X)$ of self-dual harmonic forms on $X$ when $X$ is compact (without boundary). The deformation theory of compact coassociative submanifolds with boundary and non-compact coassociative submanifolds with conical ends has been developed in \cite{KL09} and \cite{L09} respectively. 

Suppose $(M, \varphi)$ admits a fibration structure $\pi: M \to Y$ with smooth coassociative fibres. The tangent bundle $TM$ fits into an exact sequence
\[
0 \longrightarrow V \longrightarrow TM  \longrightarrow \pi^*TY  \longrightarrow 0,
\]
where $V = \ker \pi_* \subseteq TM$ is referred to as the vertical tangent bundle of the fibration. We then get a canonical decreasing filtration of the exterior tensor bundle $\Lambda^*$ of $M$ by 
\begin{equation}\label{e5.1}
\mathcal{F}^s \Lambda^tT^*_pM := \bigcap_{v_1, ..., v_s \in T_pM} \ker \left( \res \circ \iota_{v_1} \circ ... \circ \iota_{v_s} : \Lambda^{s+t}T^*_pM \longrightarrow \Lambda^tV_p^* \right),
\end{equation}
where $\res$ means the restriction map from $\Lambda^*T^*M$ to $\Lambda^*V^*$. The quotient of this filtration is identified canonically as
\begin{equation}\label{e5.2}
\mathcal{F}^{s-1}\Lambda^{s+t} / \mathcal{F}^s \Lambda^{s+t} \simeq \Lambda^s\pi^*T^*Y \otimes \Lambda^t V^*. 
\end{equation}
We denote by $\Omega^{s,t}:=\Gamma(M, \Lambda^s\pi^*T^*Y \otimes \Lambda^t V^*)$ the space of smooth sections of the corresponding quotient bundle. 

Given $y \in Y$, we write $X_y = \pi^{-1}(y)$ for the coassociative fibre over $y \in Y$ and $\vol_{X_y} = \psi|_{X_y}$ for the volume form. A $3$-form $\underline{\eta} \in \Omega^{1,2}$ of type $(1,2)$ can be regarded as a bundle map $\underline{\eta}^{\dagger}: \pi^*TY \mapsto \Lambda^2V^*$. Then the restriction of the exterior derivative of $\Omega^*(M)$ on the vertical tangent bundle $V$ gives us a vertical derivative, which we denote by $d_f: \Omega^{s,t} \to \Omega^{s, t+1}$. 

\begin{dfn}\label{d5.2}
 Let $\underline{\eta} \in \Omega^{1,2}$ be a $3$-form of type $(1,2)$. Choosing local coordinates $y = (y_1, y_2, y_3)$ on $Y$, we can write $\underline{\eta} = \sum_i \pi^*dy_i \wedge \eta_i$ with $\eta_i \in \Gamma(M, \Lambda^2V^*)$. 
\begin{enumerate}
\item We say $\underline{\eta}$ is hypersymplectic if 
\[
d_f \underline{\eta} = 0 \quad \text{ and } \quad \eta_i \wedge \eta_j|_y = 2a_{ij}(y) \vol_{X_y},
\]
where $a_{ij}(y): X_y \to \R^{3 \times 3}$ is a smooth function on $X_y$ with values in the subset of positive definite $3\times 3$ matrices. 

\item We say $\underline{\eta}$ is hyperk\"ahler if $\underline{\eta}$ is hypersymplectic and additionally $a_{ij}(y)$ is a constant map for each given $y \in Y$. 
\end{enumerate}
\end{dfn}

It is clear that the definition is independent of the choices of the local coordinates $(y_1, y_2, y_3)$. Note that a $(1,2)$-form $\underline{\eta}$ is hypersymplectic if and only if $\im \underline{\eta}^{\dagger} = \Lambda^+V^*$. In terms of local coordinates, this means that $\eta_i|_{X_y}$, $i=1,2,3$, spans the space of self-dual $2$-forms $\Lambda^+T^*X_y$ on the fibre $X_y$. When $\underline{\eta}$ is hyperk\"ahler, we get an induced metric $g_{\underline{\eta}}$ on the base $Y$ by requiring 
\begin{equation}\label{e5.3}
g_{\underline{\eta}}(u, v):=\left. \frac{\iota_{\pi^*u} \underline{\eta}\wedge \iota_{\pi^*v} \underline{\eta}}{2\vol_{X_y}}\right\vert_{x} \quad \forall u, v \in T_yY,
\end{equation}
where $x \in X_y$ is an arbitrary point. The independence of right-hand side of \eqref{e5.3} on $x$ follows from the fact that $a_{ij}(y)$ is constant for hyperk\"ahler $(1,2)$-forms in \autoref{d5.2}. 

We note that $\varphi$ defines a canonical $(1, 2)$-form $\underline{\omega} \in \Omega^{1,2}$. There are two ways to see this. One way is to identify $\varphi|_p$ with the standard positive $3$-form $\varphi_0$ on $\R^7$ and $V_p$ with the standard coassociative $\R^4$. Then it is straightforward to see that $\res \circ \iota_{v_1} \circ \iota_{v_2}(\varphi_0) = 0$ for any tangent vectors $v_1, v_2 \in \R^7$. This pointwise calculation tells us that $\varphi \in C^{\infty}(M, \mathcal{F}^0\Lambda^3)$. Then $\underline{\omega}$ is the quotient class of $\varphi$ identified in \eqref{e5.2}. Alternatively, one can choose a horizontal distribution $H$ decomposing $TM = H \oplus V$, then set $\underline{\omega}^{\dagger}(\pi^*u) = \iota_{\tilde{u}}\varphi \in C^{\infty}(M, \Lambda^2 V^*)$ for $u \in TY$ and $\tilde{u}$ the horizontal lift of $u$. The independence on the choices of the horizontal distribution $H$ follows from the fact that $\varphi|_{X_y} = 0$ for each fibre $X_y$. 

Moreover, the pointwise expression of $\varphi$ as $\varphi_0$ implies that $\im \underline{\omega}^{\dagger} = \Lambda^+V^*$ and the closedness of $\varphi$ tells us that $d_f \underline{\omega} = 0$. This shows that the induced $(1,2)$-form $\underline{\omega}$ from a closed $G_2$-structure $\varphi$ is hypersymplectic. When $\underline{\omega}$ is further assumed to be hyperk\"ahler, the discussion above tells us that $\pi:(M, g_{\varphi}) \to (Y, g_{\underline{\omega}})$ is a Riemannian submersion, which for instance justifies the statement of \cite[Proposition 3.2]{B10}. 

\begin{dfn}\label{d5.3}
A coassociative fibration $\pi: M \to Y$ of a torsion-free $G_2$-manifold $(M, \varphi)$ is said to be hyperk\"ahler if the corresponding $(1, 2)$-form $\underline{\omega}$ of $\varphi$ is hyperk\"ahler.
\end{dfn}

The metric $g_{\varphi}$ induces an Ehresmann connection (as a horizontal distribution) $H_{\varphi} := V^{\perp_{g_{\varphi}}}$.  In this way, we get a bigraded decomposition of the exterior bundle $\Lambda^*$ on $M$
\begin{equation}
\Lambda^r = \bigoplus_{s + t = r} \Lambda^{s, t}_{\varphi} := \bigoplus_{s+t=r} \Lambda^s H^*_{\varphi} \otimes \Lambda^t V^*. 
\end{equation}
Let's write $\Omega_{\varphi}^{s,t} = C^{\infty}(M, \Lambda^{s,t}_{\varphi})$ for the space of smooth sections. Then we get an orthogonal decomposition of differential forms
\[
    \Omega^r(M)=\bigoplus_{s+t=r}\Omega^{s,t}_{\varphi}. 
\]
Given a $r$-form $\omega \in \Omega^r(M)$, we write $\omega^{s,t}$ for the $(s, t)$-component of $\omega$ in $\Omega^{s,t}_{\varphi}$. 

As explained in \cite{D17}, the exterior derivative $d:\Omega^k(M)\rightarrow \Omega^{k+1}(M)$ decomposes as
\begin{equation*}
    d=d_f+d_{H_{\varphi}}+F_{H_{\varphi}}: \Omega_{\varphi}^{s,t}\longrightarrow \Omega_{\varphi}^{s,t+1}\oplus\Omega_{\varphi}^{s+1,t}\oplus\Omega_{\varphi}^{s+2,t-1}.
\end{equation*}
Here, $d_f$ is the fibrewise differential given by restricting $d$ to $V \subseteq TM$ as before; $d_{H_{\varphi}}$ is the exterior derivative on the fibration $\pi: M \to Y$ defined by coupling to the connection $H_{\varphi}$; and $F_{H_{\varphi}}$ is the curvature of $H$, which acts algebraically on forms.

\begin{dfn}
A coassociative fibration $\pi: M \to Y$ of a torsion-free $G_2$-manifold is said to be flat if the curvature of the induced Ehresmann connection vanishes, i.e. $F_{H_{\varphi}} = 0$. 
\end{dfn}

Following the notation of Donaldson \cite{D17}, we write $\underline{\omega}_{\varphi}:= \varphi^{1,2}$ and $\underline{\lambda}_{\varphi}= \varphi - \underline{\omega}_{\varphi}$. Since our choice of Ehresmann connection decomposes $TM$ orthogonally, we see that $\underline{\lambda}_{\varphi} \in \Omega^{3,0}_{\varphi}$. Since the fibres are coassociative, we know that its orthogonal complement $H_{\varphi}$ is calibrated by $\varphi$ with volume form given by $\underline{\lambda}_{\varphi}$. The closedness condition $d\varphi = 0$ is equivalent to the following equations
\begin{equation}\label{e5.5}
d_f \underline{\omega}_{\varphi} = 0 \qquad d_{H_{\varphi}} \underline{\omega}_{\varphi} = 0 \qquad d_f \underline{\lambda}_{\varphi} = -F_{H_{\varphi}} \underline{\omega}_{\varphi}. 
\end{equation}
Correspondingly, we write 
\[
\underline{\nu}_{\varphi}:= \star_{\varphi} \underline{\omega}_{\varphi} \in \Omega^{2,2}_{\varphi} \qquad \text{ and } \qquad \underline{\mu}_{\varphi}:= \star_{\varphi} \underline{\lambda}_{\varphi} \in \Omega^{0,4}_{\varphi}.
\]
Then $\psi = \underline{\mu}_{\varphi} + \underline{\nu}_{\varphi}$. We note that $\underline{\mu}_{\varphi}$ gives the volume form on fibres. The closedness condition $d\psi = 0$ is equivalent to the following equations
\begin{equation}\label{e5.6}
d_{H_{\varphi}} \underline{\mu}_{\varphi} = 0 \qquad d_{H_{\varphi}} \underline{\nu}_{\varphi} = 0 \qquad d_f \underline{\nu}_{\varphi} = -F_{H_{\varphi}} \underline{\mu}_{\varphi}. 
\end{equation}

These forms can be described locally as follows. Let $(y_1, y_2, y_3)$ be local coordinates on $U \subseteq Y$. We write $\alpha_i:= \pi^*dy_i$ for the pull-back forms on $U' = \pi^{-1}(U) \subseteq M$. Then we have the expression
\[
\underline{\omega}_{\varphi} = \sum_i \alpha_i \wedge \omega_i \qquad \underline{\lambda}_{\varphi} = \lambda \alpha_1 \wedge \alpha_2 \wedge \alpha_3,
\]
where $\omega_i \in \Omega^{0,2}(U')$ and $\lambda \in \Omega^0(U')$. Since $\underline{\nu}_{\varphi} = \star_{\varphi} \underline{\omega}_{\varphi}$, we can write $\underline{\nu}_{\varphi} = \sum_{\text{cyclic}} \alpha_i \wedge \alpha_j \wedge \nu_k$ with $\nu_i \in \Omega^{0,2}(U')$. Since each $\omega_i$ is self-dual with respect to the restricted metric on fibres, we conclude that 
\begin{equation}\label{e5.7}
\nu_i = \frac{\omega_i}{|\alpha_i|^2 \lambda}. 
\end{equation}
With these preparations, we give the following criterion for hyperk\"ahler coassociative fibrations, which is essentially demonstrated in \cite[Lemma 5]{D17}. 

\begin{lem}\label{l5.5}
Let $\pi: M \to Y$ be a coassociative fibration of a torsion-free $G_2$-manifold $(M, \varphi)$. We write $\underline{\omega} \in \Omega^{1,2}$ for the induced $(1, 2)$-form. Then the following are equivalent.
\begin{enumerate}[label=(\alph*)]
\item $\pi: M \to Y$ is hyperk\"ahler.
\item $\pi: M \to Y$ is flat. 
\item There exists a metric $g_Y$ on $Y$ so that $\pi: (M, g_{\varphi}) \to (Y, g_Y)$ is a Riemannian submersion. 
\end{enumerate}
\end{lem}

\begin{proof}
$(a) \Longrightarrow (b)$. Suppose $\pi: M \to Y$ is hyperk\"ahler. We can choose local coordinates $(y_1, y_2, y_3)$ of $U \subseteq Y$ so that $\underline{\omega}_{\varphi} = \sum_i \alpha_i \wedge \omega_i$ is hyperk\"ahler. Due to \eqref{e5.6}, it suffices to show that $d_f\underline{\nu}_{\varphi} = 0$. 

Let's write $\alpha_i^{\sharp}$ for the tangent field dual to $\alpha_i$ with respect to the metric $g_{\varphi}$. Since $\underline{\lambda}_{\varphi} = \lambda \alpha_1 \wedge \alpha_2 \wedge \alpha_3$ is the volume form on $H_{\varphi}$, we know that 
\[
\lambda = \sqrt{\det \langle \alpha^{\sharp}_i, \alpha^{\sharp}_j \rangle}
\]
Invoking \eqref{e2.22}, we see that 
\[
\left(|\alpha^{\sharp}_i|^2 \omega_i\right) \wedge \left(|\alpha^{\sharp}_j|^2\omega_j\right) \wedge \underline{\lambda}_{\varphi} 
= (\alpha^{\sharp}_i \lrcorner \; \varphi) \wedge (\alpha^{\sharp}_j\lrcorner \; \varphi) \wedge \varphi 
= 6 \langle \alpha^{\sharp}_i, \alpha^{\sharp}_j \rangle \: \underline{\lambda}_{\varphi} \wedge \underline{\mu}_{\varphi}. 
\]
Since the components in $\underline{\lambda}_{\varphi}$ are orthogonal to those in $\omega_i$ and $\underline{\mu}_{\varphi}$, we conclude that 
\begin{equation}\label{e5.8}
\omega_i \wedge \omega_j = \frac{ \langle \alpha^{\sharp}_i, \alpha^{\sharp}_j \rangle}{|\alpha^{\sharp}_i|^2 \cdot |\alpha^{\sharp}_j|^2} \: 6\underline{\mu}_{\varphi}. 
\end{equation}
By taking $i = j$, we see that $d_f|\alpha^{\sharp}_i|^2 = 0$, thus $d_f  \langle \alpha^{\sharp}_i, \alpha^{\sharp}_j \rangle = 0$. This tells us that 
\[
d_f \lambda = 0 \qquad d_f|\alpha_i|^2 = 0. 
\]
Since $d_f \omega_i = 0$ holds automatically, it follows from \eqref{e5.8} that 
\[
d\underline{\nu}_{\varphi} = \sum_{\text{cyclic}} \alpha_i \wedge \alpha_j \wedge d_f \nu_k = 0. 
\]

$(b) \Longrightarrow (a)$. When $F_{H_{\varphi}} = 0$, \eqref{e5.6} tells us that $d_f \nu_i = 0$. Combining with $d_f\lambda = 0$, \eqref{e5.8} implies that $d_f|\alpha^{\sharp}_i|^2 = 0$. Since this argument works for any local coordinates $(y_1, y_2, y_3)$ on $U$, in particular we have $d_f |\alpha^{\sharp}_i + \alpha^{\sharp}_j|^2 = 0$ by replacing $y_i$ with $y_i + y_j$. Then \eqref{e5.8} tells us that $d_f(\omega_i \wedge \omega_j / \underline{\mu}_{\varphi}) = 0$, which means $\pi: M \to Y$ is hyperk\"ahler. 

$(a) \Longrightarrow (c)$. When $\pi: M \to Y$ is hyperk\"ahler, we can choose $g_Y := g_{\underline{\omega}}$. Then the projection becomes a Riemannian submersion as discussed above \autoref{d5.3}. 

$(c) \Longrightarrow (a)$. Suppose $\pi: M \to Y$ is a Riemannian submersion. Then we can choose local coordinates $(y_1, y_2, y_3)$ around $y \in Y$ so that $(\partial_{y_1}, \partial_{y_2}, \partial_{y_3})$ forms an orthonormal frame of $T_yY$. Then $\alpha^{\sharp}_{i}|_{X_y}$ coincides with the horizontal lift of $\partial_{y_i}$. It follows from \eqref{e5.8} that $\underline{\omega}$ is hyperk\"ahler. 
\end{proof}

\autoref{l5.5} enables us to determine all hyperk\"ahler coassociative fibrations over closed $3$-manifolds. 

\begin{prop}\label{p5.6}
Let $\pi: M \to Y$ be a hyperk\"ahler coassociative fibration of a complete torsion-free $G_2$-manifold $(M, \varphi)$, where $Y$ is a closed $3$-manifold. Then 
\[
M = \R^3 \times_{G} X
\]
is the fibre product of the flat $\R^3$ and a complete hyperk\"ahler $4$-manifold $X$, where $G$ is group acting on $\R^3$ and $X$ by isometries. Moreover $Y = \R^3/G$ is a closed flat $3$-manifold. 
\end{prop}

\begin{proof}
Fix $y \in Y$. Denote by $\tilde{Y}$ the universal cover of $Y$ equipped with the pull-back metric $\tilde{g}_{\underline{\omega}}$. Using the canonical $(1,2)$-form $\underline{\omega}$, an orthonormal frame of $T_yY$ gives us a hyperk\"ahler structure on the fibre $X_y$. The flatness condition $F_{H_{\varphi}} = 0$ implies that $M$ is isomorphic to $\tilde{Y} \times_{\pi_1(Y, y)} X_y$ as a fibre bundle. Since $d_{H_{\varphi}} \underline{\omega}_{\varphi} = 0$, we further know that $\pi_1(Y, y)$ acts on $X_y$ by diffeomorphisms preserving its hyperk\"ahler structure. 

It remains to show that the metric $\tilde{g}_{\underline{\omega}}$ is flat. To see this, we choose local coordinates $(y_1, y_2, y_3)$ around $y \in Y$ so that $(\partial_{y_1}, \partial_{y_2}, \partial_{y_3})$ forms an orthonormal frame of $T_yY$. Then \eqref{e5.8} implies that $(\partial_{y_1}, \partial_{y_2}, \partial_{y_3})$ forms an orthonormal frame of $TY$ near $y$, which means $g_{\underline{\omega}}$ is flat near $y$. Alternatively, we know $\Hol(\tilde{g}_{\underline{\omega}})$ is a connected subgroup of $G_2 \cap SO(3)$, which has to be the trivial group. This also implies that $\tilde{g}_{\underline{\omega}}$ is flat. 
\end{proof}

\begin{cor}\label{c5.7}
Under the same assumption of \autoref{p5.6}, if we further assume that the isometry group of the fibre of $\pi: M \to Y$ is finite, then 
\[
M = T^3 \times_{\Gamma} X
\]
is the fibre product of a flat $3$-torus $T^3$ and a complete hyperk\"ahler $4$-manifold $X$, where $\Gamma$ is finite group acting on $T^3$ and $X$ by isometries. 
\end{cor}

\begin{proof}
Since $Y$ is a closed flat $3$-manifold, we know that $Y = T^3/ \Gamma$ for some finite group $\Gamma$ acting as isometries on $T^3$. Since the isometry group of $X$ is finite, after passing to a finite cover $\widehat{M}$, we know $\widehat{M} = T^3 \times X$. 
\end{proof}

\begin{rem}
It is unclear if one can actually find a hyperk\"ahler coassociative fibration not arising as a finite quotient of $T^3 \times X$. When $M$ is compact, this is always the case, which can be justified by the Cheeger--Gromoll's splitting theorem \cite{CG71} as well. 
\end{rem}

Let $Y = T^3$ be a flat $3$-torus, and $X$ a closed hyperk\"ahler $4$-manifold. We write $(y_1, y_2, y_3)$ for the coordinates on the universal cover $\R^3$ of $Y$. The product manifold $M = Y \times X$ admits a torsion-free $G_2$-structure 
\begin{equation}
\varphi = dy_1 \wedge dy_2 \wedge dy_3 + \sum_{i=1}^3 dy_i \wedge \omega_i.
\end{equation}
Let us write $\underline{\omega} = \sum_i dy_i \wedge \omega_i \in \Omega^{1,2}_{\varphi}$ for the hyperk\"ahler element induced by $\varphi$. Then all the main results established before applies to $\underline{\omega}$-null principal $SU(r)$-bundles.

We note that besides the $\varphi$-calibrated submanifolds $Y \times \{x\}$, there are many more $\varphi$-calibrated submanifolds in $Y \times X$ in this case. For instance, if one take a $J_i$-holomorphic curve $\Sigma$ in $X$, then $S^1 \times \Sigma \subseteq Y \times X$ is calibrated by $\varphi$.

\subsection{Cayley Fibrations}

In this subsection, we focus on $\Spin(7)$-manifolds fibred by Cayley submanifolds. This is a type of calibrated submanifold introduced in \cite{HL82}, whose local deformation theory was studied by \cite{M98}. 

We shall write $(M, \Phi)$ for a complete, torsion-free $\Spin(7)$-manifold. In particular, $M$ carries an induced metric $g_{\Phi}$ and a volume form $\vol_{\Phi}$. 

\begin{dfn}
An oriented $4$-dimensional submanifold $X \subseteq M$ is called a Cayley submanifold if $\Phi|_X = \vol_X$. 
\end{dfn}

Given a Cayley submanifold $X \subseteq M$, we identify the normal bundle $NX$ with the orthogonal complement $(TX)^{\perp_{\Phi}} \subseteq TM|_X$ with respect to the metric $g_{\Phi}$. We then have a canonical identification
\begin{equation}\label{e5.10}
\begin{split}
\Lambda^+N^*X \longrightarrow \Lambda^+T^*X. \\
\omega \longmapsto (\omega \lrcorner \; \Phi)|_X
\end{split}
\end{equation}

Now suppose $(M, \Phi)$ admits a fibration structure $\pi: M \to Y$ with smooth Cayley fibres. The metric $g_{\Phi}$ induces an Ehresmann connection $H_{\Phi} \subseteq TM$. Writing $V = \ker (\pi_*: TM \to TY)$ for the vertical sub-bundle, we get a bigraded orthogonal decomposition of differential forms as before
\begin{equation}
\Omega^r(M) = \bigoplus_{s+t = r} \Omega^{s, t}_{\Phi} =: \bigoplus_{s+t = r} C^{\infty}(M, \Lambda^s H^*_{\Phi} \otimes \Lambda^t V^*). 
\end{equation} 
The exterior derivative decomposes accordingly as
\begin{equation}
d = d_{H_{\Phi}} + d_f + F_{H_{\Phi}}: \Omega^{s,t}_{\Phi} \longrightarrow \Omega^{s+1, t}_{\Phi} \oplus \Omega^{s, t+1}_{\Phi} \oplus \Omega^{s+2, t-1}_{\Phi}. 
\end{equation}
Although the filtration \eqref{e5.1} will give us an intrinsic type decomposition, there is no canonical decomposition of $\Phi$ analogous to that of the $G_2$-calibration form $\varphi$. 

\begin{lem}
The $\Spin(7)$-structure $\Phi$ induces a canonical conformal class $c_{\Phi}$ on $Y$. 
\end{lem}

\begin{proof}
Given $y \in Y$, we choose a lift $x_o \in X_y$. Then $\pi_*: H_{\Phi, x_o} \to T_yY$ is an isomorphism. We claim that the image $\im (\pi_*: \Lambda^+H_{\Phi, x_o} \to \Lambda^2T^*Y)$ is independent of the choice of the lift $x_o$. By declaring $\Lambda^+T^*Y$ as the image, we get a conformal structure on $Y$. 

To prove the claim, we choose $\omega_o \in \Lambda^+H^*_{\Phi, x_o}$, which yields a section
\begin{equation*}
    \omega_y \in C^{\infty}(X_y, \Lambda^2 H^*_{\Phi}|_{X_y})
\end{equation*}
by lifting $\pi_*(\omega_o) \in \Lambda^2T^*Y$. This gives us a smooth function $a(x)$ on $X_y$ specified by
\[
\omega_y \wedge \omega_y = a(x) \vol_{X_y}. 
\]
Since $X_y$ is Cayler, \eqref{e5.10} tells us that $a(x)$ is positive. This proves that the images of the space of self-dual $2$-forms on $H_{\Phi}|_{X_y}$ are identified in $\Lambda^2T^*Y$. 
\end{proof}

Let us write $\underline{\omega}_{\Phi} = - \Phi^{2, 2} \in \Omega^{2, 2}_{\Phi}$ for the $(2, 2)$-component of the calibration form $\Phi$. Choosing a local frame $(\bar{\tau}_1, \bar{\tau}_2, \bar{\tau}_3)$ of $\Lambda^+T^*Y$ with respect to the conformal structure $c_{\Phi}$, we can write 
\[
\underline{\omega}_{\Phi} = \sum_{i=1}^3 \tau_i \wedge \omega_i,
\]
where $\tau_i \in \Omega^{2,0}_{\Phi}$ is the $H_{\Phi}$-lift of $\bar{\tau}_i$, and $\omega_i \in \Omega^{0,2}_{\Phi}$. 

From the standard $\Spin(7)$-form $\Phi_0$ in \eqref{e2.31}, we can decompose $\Phi$ into 
\begin{equation}
\Phi = \underline{\nu}_{\Phi} - \underline{\omega}_{\Phi} + \underline{\mu}_{\Phi} \in \Omega^{4,0}_{\Phi} \oplus \Omega^{2,2}_{\Phi} \oplus \Omega^{0,4}_{\Phi}.
\end{equation}
The self-duality of $\Phi$ implies that $\star_{\Phi} \underline{\nu}_{\Phi} = \underline{\mu}_{\Phi}$ and $\star_{\Phi} \underline{\omega}_{\Phi} = \underline{\omega}_{\Phi}$. The torsion-free condition gives us the structure equations
\begin{equation}\label{e5.14} 
d_{H_{\Phi}}\underline{\omega}_{\Phi} = 0 \qquad d_{H_{\Phi}}\underline{\mu}_{\Phi} = 0 \qquad d_f \underline{\omega}_{\Phi} = F_{H_{\Phi}}\underline{\mu}_{\Phi} \qquad d_f \underline{\nu}_{\Phi} = F_{H_{\Phi}} \underline{\omega}_{\Phi}
\end{equation}

\begin{dfn}
We say the Cayley fibration $\pi: M \to Y$ is hyperk\"ahler if near each $y \in Y$ one can choose a local frame $(\bar{\tau}_1, \bar{\tau}_2, \bar{\tau}_3)$ of $\Lambda^+T^*Y$ so that in the expression $\underline{\omega}_{\Phi} = \sum_{i=1}^3 \tau_i \wedge \omega_i$ we have 
\[
\omega_i \wedge \omega_j|_{X_y} = 2\delta_{ij} \cdot \vol_{X_y}. 
\]
\end{dfn}

We prove the analogue of \autoref{l5.5} for Cayley fibrbations. 

\begin{lem}
Let $\pi: M \to Y$ be a Cayley fibration of a complete $\Spin(7)$-manifold. Then the following are equivalent. 
\begin{enumerate}[label=(\alph*)]
\item $\pi: M \to Y$ is hyperk\"ahler.
\item $\pi: M \to Y$ is flat, i.e. $F_{H_{\Phi}} = 0$. 
\item There exists a hyperk\"ahler metric $g_Y$ on $Y$ so that $\pi: (M, g_{\varphi}) \to (Y, g_Y)$ is a Riemannian submersion. 
\end{enumerate}
\end{lem}

\begin{proof}
$(a) \Longleftrightarrow (b)$. It follows from \eqref{e5.14} that (b) is equivalent to $d_f \underline{\omega}_{\Phi} = 0$. With respect to a local frame of $\Lambda^+T^*Y$, we can write $\omega_i \wedge \omega_j|_{X_y} = 2a_{ij}(y) \cdot \vol_{X_y}$, where $(a_{ij}(y)): X_y \to \R^{3 \times 3}$ is a smooth function taking values in the space of positive definite symmetric matrices. Then $d_f \underline{\omega}_{\Phi} = 0$ is equivalent to $d a_{ij}(y) = 0$, which means we can choose a local frame to satisfy the hyperk\"ahler condition. 

$(a) \Longrightarrow (c)$. Suppose $\pi: M \to Y$ is hyperk\"ahler. Locally, we have a frame $\bar{\tau}_1, \bar{\tau}_2, \bar{\tau}_3$ of $\Lambda^+T^*Y$ and $\underline{\omega}_{\Phi} = \sum_i \tau_i \wedge \omega_i$ satisfying $\omega_i \wedge \omega_j = 2\delta_{ij} \cdot \vol_{X_y}$. Then we get a metric $g_Y$ on $Y$ characterized by 
\begin{equation}\label{e5.15}
\bar{\tau}_i \wedge \bar{\tau}_j  = 2\delta_{ij} \cdot \vol_{g_Y}. 
\end{equation}
Since any two such frames differ by an action of elements in $SO(3)$, the metric is well-defined. It follows from \eqref{e5.10} that $\tau_i \wedge \tau_j = 2 \delta_{ij} \cdot \underline{\nu}_{\Phi}$. Thus, the projection $\pi_*: H_{\Phi, x} \to T_yY$ is an isometry for each lift $x \in X_y$, which means $\pi$ is a Riemannian submersion. We note that the holonomy group $\Hol(g_Y) \subseteq \Spin(7) \cap SO(4) = SU(2)$. Thus the induced metric $g_Y$ is hyperk\"ahler. 

$(c) \Longrightarrow (a)$. Let $\bar{\tau}_1, \bar{\tau}_2, \bar{\tau}_3$ be hyperk\"ahler $2$-forms on $(Y, g_Y)$ satisfying \eqref{e5.15}. Denote by $\tau_i$ the $H_{\Phi}$-lift of $\bar{\tau}_i$, $i=1, 2, 3$. Then we can write $\underline{\omega}_{\Phi} = \sum_i \tau_i \wedge \omega_i$. Then \eqref{e2.32} implies that $\omega_i \wedge \omega_j = 2\delta_{ij} \cdot \vol_{X_y}$.
\end{proof}

The same argument as that of \autoref{p5.6} gives us the corresponding result for Cayley fibrations.

\begin{prop}\label{p5.13}
Let $\pi: M \to Y$ be hyperk\"ahler Cayley fibration of a complete torsion-free $\Spin(7)$-manifold $(M, \Phi)$. Then we get a hyperk\"ahler metric $g_Y$ and a hyperk\"ahler $4$-manifold $(X, g_X)$ so that $(M, g_{\Phi})$ is isometric to $\tilde{Y} \times_{\pi_1(Y)} X$, where $\tilde{Y}$ is the universal cover of $Y$ and $\pi_1(Y)$ acts on $X$ via diffeomorphisms preserving its hyperk\"ahler structure.  
\end{prop}

Let $(Y, \tau_1, \tau_2, \tau_3)$ and $(X, \omega_1, \omega_2, \omega_3)$ be two compact hyperk\"ahler $4$-manifolds. Their product $M = Y \times X$ is equipped with a torsion-free $\Spin(7)$-structure
\begin{equation}
\Phi:= \vol_Y + \vol_X - \sum_{i=1}^3 \tau_i \wedge \omega_i. 
\end{equation}
Then the projection $\pi_1: M \to Y$ is a hyperk\"ahler Cayley fibration. Let us write $\underline{\eta} = \vol_X - \sum_{i=1}^3 \tau_i \wedge \omega_i$. Then all the main results proved in \autoref{s4} applies in this case by setting $\alpha = \underline{\eta}$. 

\section{More Examples}\label{s6}

In this section, we discuss the remaining examples that appeared in \autoref{t1}.

\subsection{$G_2$-Instantons over $S^1 \times Z$}

We start with a codimension-$1$ case of dimension reduction, which recovers and extends the results of Wang \cite{W20} in a more explicit form. 
 
 \begin{dfn}
Let $Z$ be a smooth $6$-manifold. A half-closed $SU(3)$-structure on $Z$ is a quadruple $(J, g, \omega, \Omega)$ consisting of 
\begin{enumerate}
\item an almost complex structure $J$, 
\item an Hermitian metric $g$ with respect to $J$,
\item a positive real $(1, 1)$-form $\omega$ given by $\omega(-, -):= g(J-, -)$,
\item a nowhere vanishing $(3, 0)$-form $\Omega$,
\end{enumerate}
which satisfies the conditions 
\[
\frac{\omega^3}{6} = \frac{1}{4} \Rea \Omega \wedge \Ima \Omega \qquad \text{ and } \qquad d(\Rea \Omega) = 0.
\]
\end{dfn}

To simplify notation, we shall write $(\omega, \Omega)$ for an $SU(3)$-structure and omit $J$ and $g$. Let $Z$ be a closed smooth $6$-manifold equipped with a half-closed $SU(3)$-structure $(\omega, \Omega)$. Then we get a closed $7$-manifold $M = S^1 \times Z$ equipped with a closed $G_2$-structure
\begin{equation}\label{e6.1}
\varphi = dt \wedge \omega + \Rea \Omega,
\end{equation}
where $t$ is the arc-length parameter on the circle $S^1$ of length $2\pi$. For simplicity, we shall omit writing out the pull-back maps $\pi_1^*$ and $\pi_2^*$ in the notation for differential forms later when there is no ambiguity. 

Consider a principal $U(r)$-bundle $\pi: P \to M$. We fix a reference point $t_0 \in S^1$ and identify $P|_{Z_{t_0}} \simeq P_Z$ with a $U(r)$-bundle over $Z$, where $Z_t := \{t\} \times Z$ denotes the fibre over $t \in S^1$. Recall that the degree of $P_Z$ is defined to be $\deg P_Z = \langle c_1(P_Z) \smile [\omega]^2, [Z] \rangle$. We start with a simple observation that is covered in \cite[Theorem 1.17]{W20} when $P = \pi_2^*P_Z$ is the pull-back bundle. 

\begin{lem}\label{l6.2} 
Let $M = S^1 \times Z$ be the product $G_2$-manifold as above, and $\pi: P \to M$ a $U(r)$-bundle (not necessarily $\Rea \Omega$-null). Suppose the moduli space of $G_2$-instantons $\M_{\varphi}(M, P)$ is non-empty. Then $\deg P_Z = 0$. 
\end{lem}

\begin{proof}
Let $A$ be a $G_2$-instanton on $P$. Then \eqref{e2.30} gives us 
\begin{equation}
\star_{\varphi}(F_A \wedge \star_{\varphi} \varphi) = F_A \;\lrcorner\; \varphi = 0. 
\end{equation}
We pick $I_{\epsilon} = (t_0 - \epsilon, t_0 + \epsilon) \subseteq S^1$ a small neighborhood of $t_0$ in $S^1$, and choose an identification $P|_{I_{\epsilon} \times Z} \simeq I_{\epsilon} \times P_Z$. Then over $I_{\epsilon} \times Z$, we can write $A = B_t + \sigma dt$, where $B_t$ is a family of connection on $P_Z$, and $\sigma(t, -) \in \Omega^0(Z, \g_{P_Z})$. The curvature of $F_A$ then takes the form 
\begin{equation}
F_A = F_{B_t} + dt \wedge \left(d_{B_t} \sigma - \dot{B}_t \right). 
\end{equation}
Since $\varphi = dt \wedge \omega + \Rea \Omega$, the $(1, 0)$-component of $F_A \;\lrcorner\; \varphi$ is $ dt \wedge (F_{B_t} \;\lrcorner\; \omega)$, which must vanish for each $t \in I_{\epsilon}$. Thus, we conclude that 
\begin{equation}
0 = \int_{Z_{t_0}} \tr (F_{B_{t_0}}) \;\lrcorner\; \omega \vol_{Z_{t_0}}  = \frac{1}{2} \int_{Z_t} \tr (F_{B_{t_0}}) \wedge \omega^2 = -\pi i \langle c_1(P_Z) \smile [\omega]^2, [Z] \rangle, 
\end{equation}
which finishes the proof. 
\end{proof}

We summarize our main results in this case as follows. 

\begin{thm}\label{t6.3}
Let $M =S^1 \times Z$ be the product $G_2$-manifold of a circle and a closed $6$-manifold $Z$ equipped with a half-closed $SU(3)$-structure $(\omega, \Omega)$. Suppose $\pi: P \to M$ is a $U(r)$-bundle satisfying $\deg P|_{Z_{t_0}} = 0$ for a fixed basepoint $t_0 \in S^1$. 
\begin{enumerate}[label=(\alph*)]
\item If $P$ is $\Rea \Omega$-negative, then $\M_{\varphi}(M, P) = \varnothing$. 
\item If $P$ is $\Rea \Omega$-null, then each $G_2$-instanton $A$ on $P$ satisfies 
\[
F_A^{2, 0} = 0 \quad \text{ and } \quad F_A^{1,1} = 0. 
\]
\end{enumerate}
Moreover, the moduli space $\M_{\varphi}(M, P)$ is non-empty if and only if
\begin{enumerate}
\item[(b.1)] The moduli space $\M_{\omega}(Z, P_Z)$ of HYM connections on $P_Z$ is non-empty.
\item[(b.2)] The family bundle is $\underline{P} = \R \times_{\rho} P_Z$ for some representation $\rho:\pi_1(S^1, t_0) \to \G_{P_Z}$ satisfying $\im \rho \subseteq \Stab(B)$ for some HYM connection $B$ on $P_Z$. 
\end{enumerate}
\end{thm}

The pull-back bundle $\pi_2^*P_Z \to M$ considered in \cite{W20} satisfies the $\Rea \Omega$-null assumption automatically, which corresponds to the trivial representation $\rho: \pi_1(S^1, t_0) \to \G_{P_Z}$. When the structure group is $U(r)$, irreducible connections are defined to be those having $U(1)$ as stabilizers. Since the Picard torus $\Pic(S^1)$ can be identified with the circle $S^1$ itself, the twisting map \eqref{e4.17} is defined on $\Pic(S^1) \times \mathcal{B}^*_{P_Z}$ and restricts as a homeomorphism between $S^1 \times \M_{\omega}(Z, P_Z)$ and $\M_{\varphi}(M, P)$, following the first part proof of \autoref{t1.1}. Thus, the main result of Wang \cite[Theorem 1.17 $\mathbb{I}$]{W20} can be recovered from \autoref{t6.3}.

To get an elliptic deformation theory for $G_2$-instantons, the calibration form $\varphi$ is required to be coclosed. To this end, we consider Calabi--Yau fibres. Note that a half-closed $SU(3)$-structure $(J, g, \omega, \Omega)$ defines a Calabi--Yau structure if we further require 
\begin{equation}
(a) \; J \text{ is integrable} \quad  (b) \;  \omega \text{ is closed} \quad (c) \; \Omega \text{ is holomorphic}
\end{equation}
In particular, the $G_2$-structure $\varphi$ in \eqref{e6.1} is torsion-free when $Z$ is a Calabi--Yau $3$-fold. Then the main results hold for $\Rea\Omega$-null bundles. 

\subsection{$\Spin(7)$-Instantons over $S^1 \times X$}

Let $(X, \varphi)$ be a closed, torsion-free $G_2$-manifold. The product manifold $M = S^1 \times X$ is equipped with a $\Spin(7)$-structure 
\begin{equation}
\Phi = dt \wedge \varphi + \psi.
\end{equation}
As explained in \autoref{s2}, we know that $\lambda_{\Phi} = \lambda_{\varphi} = -1$. Thus the our main results applied to $\varphi$-null bundles covers the second part of Wang’s main result \cite[Theorem 1.17 $\mathbb{II}$]{W20} for the same reason explained in the paragraph following \autoref{t6.3}.

\subsection{Product K\"ahler Manifolds}

We investigate the last case in \autoref{t1}. Let $(Y, \omega_1)$ and $(X, \omega_2)$ be compact K\"ahler manifolds of complex dimension $m_1$ and $m_2$, respectively. Let $M = Y \times X$ be the product K\"ahler manifold equipped with the K\"ahler form $\omega = \pi_1^* \omega_1 + \pi_2^* \omega_2$. Let $\pi: P \to M$ be a principal $SU(r)$-bundle whose restriction on the fibres is identified with a fixed bundle $\pi_X: P_X \to X$.

To comply with the set-up of \autoref{s4}, we let
\begin{equation}
\theta = \omega^{m-2} = \vol_Y \wedge \frac{(m-2)!}{(m_2-2)!} \omega_2^{m_2 - 2} + \alpha = \vol_Y \wedge \beta + \alpha,
\end{equation}
where $m = m_1 + m_2$. It follows from \eqref{e2.21} that the absolute minimum of the Yang--Mills functional on $\A_P$ is 
\begin{equation}\label{e6.8}
\frac{16r\pi^2}{(m-2)!} \langle \kappa(P) \smile [\omega]^{m-2}, [M] \rangle = \frac{16r\pi^2}{(m-2)!}  \langle \kappa(P) \smile [\theta], [M] \rangle,
\end{equation}
while the absolute minimum of the Yang--Mills functional on $\A_{P_X}$ is 
\begin{equation}\label{e6.9}
\frac{16r\pi^2}{(m_2-2)!} \langle \kappa(P_X) \smile [\omega_2]^{m_2-2}, [X] \rangle =  \frac{16r\pi^2}{(m-2)!} \langle \kappa(P_X) \smile [\beta], [X] \rangle.
\end{equation}
Comparing with \autoref{p2.3}, we see that $\lambda_{\theta} = \eta_{\beta} = -1$. Thus all our main results applies to this case as well.

\bibliographystyle{plain}
\bibliography{RefG2}

\begin{thebibliography}{10}

\bibitem{A88}
Michael Atiyah.
\newblock Topological quantum field theories.
\newblock {\em Inst. Hautes \'{E}tudes Sci. Publ. Math.}, (68):175--186, 1988.

\bibitem{B01}
Scott Baldridge.
\newblock Seiberg-{W}itten invariants of 4-manifolds with free circle actions.
\newblock {\em Commun. Contemp. Math.}, 3(3):341--353, 2001.

\bibitem{B10}
D.~Baraglia.
\newblock Moduli of coassociative submanifolds and semi-flat {$G_2$}-manifolds.
\newblock {\em J. Geom. Phys.}, 60(12):1903--1918, 2010.

\bibitem{CG71}
Jeff Cheeger and Detlef Gromoll.
\newblock The splitting theorem for manifolds of nonnegative {R}icci curvature.
\newblock {\em J. Differential Geometry}, 6:119--128, 1971/72.

\bibitem{DK90}
S.~K. Donaldson and P.~B. Kronheimer.
\newblock {\em The geometry of four-manifolds}.
\newblock Oxford Mathematical Monographs. The Clarendon Press, Oxford
  University Press, New York, 1990.
\newblock Oxford Science Publications.

\bibitem{DT98}
S.~K. Donaldson and R.~P. Thomas.
\newblock Gauge theory in higher dimensions.
\newblock In {\em The geometric universe ({O}xford, 1996)}, pages 31--47.
  Oxford Univ. Press, Oxford, 1998.

\bibitem{D17}
Simon Donaldson.
\newblock Adiabatic limits of co-associative {K}ovalev-{L}efschetz fibrations.
\newblock In {\em Algebra, geometry, and physics in the 21st century}, volume
  324 of {\em Progr. Math.}, pages 1--29. Birkh\"{a}user/Springer, Cham, 2017.

\bibitem{DS11}
Simon Donaldson and Ed~Segal.
\newblock Gauge theory in higher dimensions, {II}.
\newblock In {\em Surveys in differential geometry. {V}olume {XVI}. {G}eometry
  of special holonomy and related topics}, volume~16 of {\em Surv. Differ.
  Geom.}, pages 1--41. Int. Press, Somerville, MA, 2011.

\bibitem{GM24}
D.~Galt and L.~Ma.
\newblock Adiabatic limit and instantons over product manifolds.
\newblock {\em In preparation}, 2025.

\bibitem{HL82}
Reese Harvey and H.~Blaine Lawson, Jr.
\newblock Calibrated geometries.
\newblock {\em Acta Math.}, 148:47--157, 1982.

\bibitem{IN90}
Mitsuhiro Itoh and Hiraku Nakajima.
\newblock Yang-{M}ills connections and {E}instein-{H}ermitian metrics.
\newblock In {\em K\"{a}hler metric and moduli spaces}, volume 18-{\rm II} of
  {\em Adv. Stud. Pure Math.}, pages 395--457. Academic Press, Boston, MA,
  1990.

\bibitem{J96}
Dominic~D. Joyce.
\newblock Compact {R}iemannian {$7$}-manifolds with holonomy {$G_2$}. {I},
  {II}.
\newblock {\em J. Differential Geom.}, 43(2):291--328, 329--375, 1996.

\bibitem{J00}
Dominic~D. Joyce.
\newblock {\em Compact manifolds with special holonomy}.
\newblock Oxford Mathematical Monographs. Oxford University Press, Oxford,
  2000.

\bibitem{KL21}
Spiro Karigiannis and Jason~D. Lotay.
\newblock Bryant-{S}alamon {$\rm G_2$} manifolds and coassociative fibrations.
\newblock {\em J. Geom. Phys.}, 162:Paper No. 104074, 60, 2021.

\bibitem{KL09}
Alexei Kovalev and Jason~D. Lotay.
\newblock Deformations of compact coassociative 4-folds with boundary.
\newblock {\em J. Geom. Phys.}, 59(1):63--73, 2009.

\bibitem{LM89}
H.~Blaine Lawson, Jr. and Marie-Louise Michelsohn.
\newblock {\em Spin geometry}, volume~38 of {\em Princeton Mathematical
  Series}.
\newblock Princeton University Press, Princeton, NJ, 1989.

\bibitem{L09}
Jason~D. Lotay.
\newblock Deformation theory of asymptotically conical coassociative 4-folds.
\newblock {\em Proc. Lond. Math. Soc. (3)}, 99(2):386--424, 2009.

\bibitem{M98}
Robert~C. McLean.
\newblock Deformations of calibrated submanifolds.
\newblock {\em Comm. Anal. Geom.}, 6(4):705--747, 1998.

\bibitem{MOY97}
Tomasz Mrowka, Peter Ozsv\'{a}th, and Baozhen Yu.
\newblock Seiberg-{W}itten monopoles on {S}eifert fibered spaces.
\newblock {\em Comm. Anal. Geom.}, 5(4):685--791, 1997.

\bibitem{M99}
Vicente Mu\~{n}oz.
\newblock Ring structure of the {F}loer cohomology of {$\Sigma\times{\bf
  S}^1$}.
\newblock {\em Topology}, 38(3):517--528, 1999.

\bibitem{S14}
Henrique~N. S\'{a}~Earp.
\newblock Generalised {C}hern-{S}imons theory and {$\rm G_2$}-instantons over
  associative fibrations.
\newblock {\em SIGMA Symmetry Integrability Geom. Methods Appl.}, 10:Paper 083,
  11, 2014.

\bibitem{S89}
Simon Salamon.
\newblock {\em Riemannian geometry and holonomy groups}, volume 201 of {\em
  Pitman Research Notes in Mathematics Series}.
\newblock Longman Scientific \& Technical, Harlow; copublished in the United
  States with John Wiley \& Sons, Inc., New York, 1989.

\bibitem{TT04}
Terence Tao and Gang Tian.
\newblock A singularity removal theorem for {Y}ang-{M}ills fields in higher
  dimensions.
\newblock {\em J. Amer. Math. Soc.}, 17(3):557--593, 2004.

\bibitem{T82}
Clifford~Henry Taubes.
\newblock Self-dual {Y}ang-{M}ills connections on non-self-dual
  {$4$}-manifolds.
\newblock {\em J. Differential Geometry}, 17(1):139--170, 1982.

\bibitem{T11}
Clifford~Henry Taubes.
\newblock {\em Differential geometry}, volume~23 of {\em Oxford Graduate Texts
  in Mathematics}.
\newblock Oxford University Press, Oxford, 2011.
\newblock Bundles, connections, metrics and curvature.

\bibitem{T00}
Gang Tian.
\newblock Gauge theory and calibrated geometry. {I}.
\newblock {\em Ann. of Math. (2)}, 151(1):193--268, 2000.

\bibitem{U82}
Karen~K. Uhlenbeck.
\newblock Connections with {$L\sp{p}$} bounds on curvature.
\newblock {\em Comm. Math. Phys.}, 83(1):31--42, 1982.

\bibitem{W13b}
Thomas Walpuski.
\newblock {$\rm G_2$}-instantons on generalised {K}ummer constructions.
\newblock {\em Geom. Topol.}, 17(4):2345--2388, 2013.

\bibitem{W17}
Thomas Walpuski.
\newblock {$G_2$}-instantons, associative submanifolds and {F}ueter sections.
\newblock {\em Comm. Anal. Geom.}, 25(4):847--893, 2017.

\bibitem{W17a}
Thomas Walpuski.
\newblock {$\rm Spin(7)$}-instantons, {C}ayley submanifolds and {F}ueter
  sections.
\newblock {\em Comm. Math. Phys.}, 352(1):1--36, 2017.

\bibitem{W22}
Donghao Wang.
\newblock On finite energy monopoles on {$\Bbb C \times\Sigma$}.
\newblock {\em Comm. Anal. Geom.}, 30(2):381--449, 2022.

\bibitem{W20}
Yuanqi Wang.
\newblock Moduli spaces of {$G_2$}-instantons and {$Spin(7)$}-instantons on
  product manifolds.
\newblock {\em Ann. Henri Poincar\'{e}}, 21(9):2997--3033, 2020.

\end{thebibliography}

\end{document}